\documentclass[12pt,a4paper,oneside,reqno,notitlepage]{amsart}
\usepackage{amssymb}

\usepackage{amsfonts}
\usepackage{amscd,amsmath}
\usepackage{euscript}
\usepackage{comment}
\usepackage[PostScript=dvips]{diagrams}

\usepackage[arrow,matrix,curve]{xy}\SilentMatrices
\def\xyma{\xymatrix@M.7em}

\numberwithin{equation}{section}

\usepackage[T2A]{fontenc}
\newtheorem{cor}{Corollary}[section]

\newtheorem{prop}{Proposition}[section]
\newtheorem{theorem}{Theorem}[section]
\newtheorem{lemma}{Lemma}[section]
\newtheorem{remark}{Remark}[section]
\newtheorem{example}{Example}[section]

\def\dim{\mathrm{dim}}

\def\Z{{\mathbb{Z}}}

\def\Len{\mathrm{Len}}

\def\bee{\begin{equation}}
\def\ee{\end{equation}}

\def\Tor{\mathrm{Tor}}

\def\RP{\mathbb{R}\mathrm{P}}
\def\CP{\mathbb{C}\mathrm{P}}
\def\sk{\mathrm{sk}}
\def\RP{\mathbb{R}\mathrm{P}}

\begin{document}
\title{On homotopy groups of the suspended classifying spaces}
\author{Roman Mikhailov}

\address{Steklov Mathematical Institute, Gubkina 8, Moscow, Russia 119991}
\email{romanvm@mi.ras.ru}

\author{Jie Wu}

\address{Department of Mathematics, National University of Singapore, 2 Science Drive 2
Singapore 117542} \email{matwuj@nus.edu.sg}
\urladdr{www.math.nus.edu.sg/\~{}matwujie}
\thanks{The research of the first author is partially supported by RFBR (grant 08-01-91300 $IND_a$). The research of the second author is partially supported by the Academic Research Fund
of the National University of Singapore R-146-000-101-112}

\begin{abstract}
In this paper, we determine the homotopy groups $\pi_4(\Sigma
K(A,1))$ and $\pi_5(\Sigma K(A,1))$ for abelian groups $A$ by using different facts and methods from group
theory and homotopy theory: derived functors, the Carlsson simplicial
construction, the Baues-Goerss spectral sequence, homotopy
decompositions and the methods of algebraic K-theory. As the applications, we also determine $\pi_i(\Sigma K(G,1))$ with $i=4,5$ for some non-abelian groups $G=\Sigma_3$ and $\mathrm{SL}(\Z)$, and $\pi_4(\Sigma K(A_4,1))$ for the $4$-th alternating group $A_4$.
\end{abstract}

\maketitle

\section{Introduction}

It is well-known that the suspension functor applied to a topological
space shifts homology groups, but "chaotically" changes homotopy
groups. For example, one can take a circle $S^1$, whose homotopy
type is very simple. Its suspension $\Sigma S^1=S^2$ has obvious
homology groups, however the problem of investigating the homotopy
groups of $S^2$ is one of the deepest problems of algebraic
topology. Consider the following functors from the category of
groups to the category of abelian groups:
$$
\pi_n(\Sigma^mK(-,1)):{\sf Gr\to Ab},\ n\geq 1, m\geq 1
$$
defined by $
A\mapsto \pi_n(\Sigma^mK(A,1)),
$
where $\Sigma^m$ is the $m$-fold suspension. It is clear that $
\pi_n(\Sigma^mK(\mathbb Z,1))=\pi_n(S^{m+1}),$ that is the homotopy
groups of spheres appear as the simplest case of a general theory
of homotopy groups of suspensions of classifying spaces.

For the case $m=1,2$ and $n=3,4$ there is the following natural
commutative diagram with exact rows \cite{BL}:
\begin{equation}\label{bldia}\xyma{ 0\ar@{->}[r] & \pi_3(\Sigma K(G,1))\ar@{->}[r]
\ar@{->>}[d]& G\otimes G
\ar@{->}[r] \ar@{->>}[d] & [G,G]\ar@{->}[r] \ar@{=}[d] & 1\\
0 \ar@{->}[r] & \pi_4(\Sigma^2K(G,1)) \ar@{->>}[d] \ar@{->}[r] &
G\widetilde \otimes
G \ar@{->}[r] \ar@{->>}[d] & [G,G]\ar@{->}[r] \ar@{=}[d] & 1\\
0 \ar@{->}[r] & H_2(G) \ar@{->}[r] & G\wedge G \ar@{->}[r] &
[G,G]\ar@{->}[r] & 1}
\end{equation}
where $G\otimes G$ is the non-abelian square of $G$ in the sense of Brown-Loday~\cite{BL}, $G\widetilde
\otimes G$ (resp. $G\wedge G$) is the quotient of $G\otimes G$ by
the normal subgroup generated by elements $g\otimes h+h\otimes g$
(resp. $g\otimes g, g\in G$). In particular, for an abelian group
$A$, there are natural isomorphisms
\begin{align*}
& \pi_3(\Sigma K(A,1))\simeq A\otimes A\\
& \pi_4(\Sigma^2K(A,1))\simeq \pi_2^SK(A,1)\simeq A\widetilde
\otimes A.
\end{align*}

The purpose of this article is to determine the homotopy groups $\pi_4(\Sigma K(A,1))$ and
$\pi_5(\Sigma K(A,1))$ for abelian groups $A$. In order to investigate the structure of these
homotopy groups, we use different facts and methods of group
theory and homotopy theory: derived functors, the Carlsson simplicial
construction, the Baues-Goerss spectral sequence~\cite{BG}, homotopy
decompositions and the methods of algebraic K-theory. The
combination of these different methods provides an effective way
for determining these homotopy groups. As reader will see, some
our computations use commutator tricks in simplicial groups.

The homotopy group $\pi_4(\Sigma K(A,1))$ as a functor on $A$ can be given as follows:

\begin{theorem}[Theorem~\ref{theorem3.1}]\label{theorem1.1}
Let $A$ be any abelian group. Then there is a natural short exact
sequence
$$
(\Lambda^2(A)\otimes A)^{\oplus 2}\oplus A\otimes A\otimes
\mathbb Z/2\rInto \pi_4(\Sigma K(A,1))\rOnto \Tor(A,A).
$$
Moreover  $(\Lambda^2(A)\otimes A)^{\oplus 2}$ is an (unnatural) summand of $\pi_4(\Sigma K(A,1))$.
\end{theorem}

An interesting point of this theorem is that the functor $\pi_4(\Sigma K(A,1))$ has $\Tor(A,A)$ as a natural quotient. For determining the structure of the group $\pi_4(\Sigma K(A,1))$, one has to solve the group extension problem in Theorem~\ref{theorem1.1}. For finitely generated abelian groups $A$, we are able to solve this problem. Given a finitely generated abelian group $A$, let
$$
A=A_1\oplus\bigoplus_{
\begin{array}{c}
r\geq 1\\
p\textrm{ is a prime}\\
\end{array}}A_{p^r}
$$
be the primary decomposition of $A$, where $A_1$ is torsion free and $A_{p^r}$ is a free $\Z/p^r$-module.

\begin{theorem}[Theorem~\ref{theorem3.2}]\label{theorem1.2}
Let $A$ be any finitely generated abelian group. Let $A=A_2\oplus B$ with $B=A_1\oplus \bigoplus_{p^r\not=2}A_{p^r}$. Then
$$
\begin{array}{rcl}
\pi_4(\Sigma K(A,1))&\cong& \frac{1}{2}(A_2\otimes A_2)\oplus (A_2\otimes B)^{\oplus 2}\oplus B^{\otimes 2}\otimes\Z/2\oplus (A\otimes \Lambda^2(A))^{\oplus 2}\\
 & & \ \ \oplus \Tor(A_2,B)^{\oplus 2}\oplus \Tor(B,B),\\
\end{array}
$$
where $\frac{1}{2}(A_2\otimes A_2)$ is a free $\Z/4$-module with rank of $\dim_{\Z/2}(A_2\otimes A_2)$.
\end{theorem}

One point of this theorem is that the (maximal) elementary $2$-group summand $A_2$ of $A$ plays a key role in the group extension problem. Roughly speaking $A_2\otimes A_2$ is half down in the group $\pi_4(\Sigma K(A,1))$.

As the applications of Theorems~\ref{theorem1.1} and~\ref{theorem1.2}, we are able to compute $\pi_4(M(\Z/2^r,2))$ and their connections with $\pi_4(\Sigma K(\Z/2^r,1))$. As the direct consequences, the homotopy groups $\pi_4(\Sigma \RP^n)$  and $\pi_4(\Sigma K(\Sigma_3,1))$ are determined. (See subsection 3.2 for the computations of these homotopy groups.)

For the homotopy group $\pi_5(\Sigma K(A,1))$, as a functor, it can be described by two exact sequences given in diagram~(\ref{e433}). Unfortunately it seems too complicated to produce a canonical functorial short exact sequence description for the functor $\pi_5(\Sigma K(A,1))$ from diagram~(\ref{e433}). For any finitely generated abelian group $A$, we determine $\pi_5(\Sigma K(A,1))$ in an un-functorial way by the following steps:
\begin{enumerate}
\item[1)] From the Hopf fibration, $\pi_5(\Sigma K(A,1))\cong \pi_5(\Sigma K(A,1)\wedge K(A,1))$;
\item[2)] Take a primary decomposition of $A$ and write $K(A,1)$ as a product of copies of $S^1=K(\Z,1)$ and $K(\Z/p^r,1)$;
\item[3)] By using the fact that
$$
\Sigma X\times Y\simeq \Sigma X\vee \Sigma Y\vee \Sigma X\wedge Y,
$$
write $\Sigma K(A,1)\wedge K(A,1)$ as a wedge of the spaces in the form
$$
X=\Sigma^m K(\Z/p_1^{r_1},1)\wedge K(\Z/p_2^{r_2},1)\wedge\cdots \wedge K(\Z/p_t^{r_t},1)
$$

with $m+t\geq 3$ and $m\geq 1$; \item[4)] By applying the
Hilton-Milnor Theorem, $\pi_5(\Sigma K(A,1))$ becomes a summation
of $\pi_5(X)$ for some $X$ in the above form.
\end{enumerate}
For the spaces $X$ in the above form, it is contractible if $p_i\not=p_j$ for some $i\not=j$ and $\pi_5(X)$ can be determined in Proposition 4.2 for an odd prime $p$. The only difficult part is to compute $\pi_5(X)$ for $X$ given in the form
$$
X=\Sigma^m K(\Z/2^{r_1},1)\wedge \cdots \wedge K(\Z/2^{r_t},1)
$$
with $m+t\geq 3$ and $m\geq 1$. Our computations are then given case-by-case (Propositions~\ref{proposition4.4}-\ref{proposition4.6} and Theorems~\ref{theorem4.1}-\ref{theorem4.4}), in which different methods are involved. An instructional example is as follows:

Let $A=\Z\oplus \Z/2$. According to 1), $\pi_5(\Sigma K(A,1))\cong \pi_5(\Sigma K(A,1)\wedge K(A,1))$. As in 3),
$$
\begin{array}{rcl}
\Sigma K(A,1)\wedge K(A,1)&\simeq& \Sigma (S^1\times \RP^\infty)\wedge (S^1\times \RP^\infty)\\
&\simeq& \Sigma (S^1\vee \RP^\infty\vee \Sigma \RP^\infty)\wedge (S^1\vee \RP^\infty\vee \Sigma\RP^\infty)\\
&=&S^3\vee\bigvee\limits^2\Sigma^2\RP^\infty \vee\bigvee\limits^2\Sigma^3\RP^\infty\vee \Sigma \RP^\infty\wedge\RP^\infty\vee\\
& & \ \ \vee\bigvee\limits^2\Sigma^2\RP^\infty\wedge\RP^\infty \vee\Sigma^3\RP^\infty\wedge\RP^\infty.\\
\end{array}
$$
By applying the Hilton-Milnor Theorem as in 4), $\pi_5(\Sigma K(A,1)\wedge K(A,1))$ is a summation of
$$
\pi_5(S^3), \ \pi_5(\Sigma^2\RP^\infty), \pi_5(\Sigma \RP^\infty\wedge \RP^\infty),\pi_5(\Sigma (\RP^\infty)^{\wedge 3}),\cdots
$$
with multiplicities. From Theorem~\ref{theorem4.4}, we have
$$
\pi_5(\Sigma^2\RP^\infty)=\Z/8
$$
and by Proposition~\ref{proposition4.6} and Theorem~\ref{theorem4.1}, we have
$$
\pi_5(\Sigma \RP^\infty\wedge \RP^\infty)=\pi_5(\Sigma (\RP^\infty)^{\wedge 3})=\Z/2^{\oplus 2}.
$$
The group $\pi_5(\Sigma K(\Z\oplus\Z/2))$ will be determined by filling all possible summands with multiplicities.

As the applications of our computations on $\pi_5(\Sigma K(A,1))$, we are able to determine
$\pi_5(\Sigma \RP^n)$ (Proposition~\ref{proposition4.8}) and $\pi_5(\Sigma K(\Sigma_3,1))$ (Proposition~\ref{proposition4.9}).

In section 2 we recall certain facts from the homotopy theory, such as the Whitehead exact sequence, the Carlsson simplicial construction
and describe a spectral sequence (\ref{specseq}), which converges
to $\pi_*(\Sigma^mK(A,1))$ for any abelian group $A$, with $E^2$-terms are given by the derived functors of certain
polynomial functors. We illustrate how it works in Theorem 4.4 for computing
$$
\pi_5(\Sigma^2 K(\Z/2^r,1))=
\begin{array}{lcl}
\mathbb Z/8&\textrm{ if }&r=1,\\
\Z/2^{r+1}\oplus\Z/2&\textrm{ if}&r>1.\\
\end{array}
$$
The interesting point is of course how $\Z/8$ shows up in the case $r=1$ while it becomes $\Z/2^{r+1}\oplus \Z/2$ for $r>1$. The proof is also based on the computations of the derived functors of the antisymmetric
square $\widetilde\otimes^2$.

There is a natural relation between the problem considered and
algebraic K-theory. Since the plus-construction $K(G,1)\to
K(G,1)^+$ is a homological equivalence, there is a natural weak
homotopy equivalence
$$
\Sigma K(G,1)\to \Sigma (K(G,1)^+)
$$
This defines the natural suspension map:
$$
\pi_n(K(G,1)^+)\to \pi_{n+1}(\Sigma(K(G,1)^+))=\pi_{n+1}(\Sigma
K(G,1))
$$
for $n\geq 1$.
This map was studied in \cite{BC} in the case of a perfect group
$G$. We consider the case $G=E(R)$, i.e. the group of elementary
matrices over a ring $R$. In this case the natural map
$$
K_3(R)=\pi_3(K(E(R),1)^+)\to \pi_4(\Sigma K(E(R),1))
$$
is an isomorphism (Theorem 5.1). The natural relation to
K-theory gives a way how to compute homotopy groups $\pi_i(\Sigma
K(E(R),1))$ for $i=4,5$ for some rings. For example, the case
$G=SL(\mathbb Z)$ is considered. As an application of our methods, we also determine that
$\pi_4(\Sigma K(A_4,1))=\Z/4$ for the $4$-th alternating group $A_4$.

The article is organized as follows. We give a brief review for the quadratic functors and the simplicial resolutions in Section 2. The determination of $\pi_4(\Sigma K(A,1))$ is given in section 3, where the proofs of Theorems~\ref{theorem1.1} and ~\ref{theorem1.2} are also given. In section 4, we give case-by-case computations for $\pi_5(\Sigma K(A,1))$. In Section 5, we give some relations to $K$-theory.

\section{The Quadratic Functors and the Simplicial Resolutions}
\subsection{Whitehead Quadratic Functor}
In~\cite[Chapter II]{Wh}, J. H. C Whitehead introduce the \textit{universal quadratic functor} $\Gamma_2$ from abelian groups to abelian groups as follows: Let $A$ be any abelian group. Then $\Gamma_2(A)$ is the group generated by the symbols $\gamma(x)$, one for each $x\in A$, subject to the defining relations
\begin{enumerate}
\item $\gamma(-x)=\gamma(x)$;
\item $\gamma(x+y+z)-\gamma(x+y)-\gamma(y+z)-\gamma(x+z)+\gamma(x)+\gamma(y)+\gamma(z)=0.$
\end{enumerate}
\textbf{Note.} According to~\cite[p. 61]{Wh}, the group $\Gamma_2(A)$ is abelian and so the multiplication in $\Gamma_2(A)$ is denoted by $+$. Define
$$
\gamma(x,y)=\gamma(x+y)-\gamma(x)-\gamma(y).
$$
The following proposition helps for determining the group $\Gamma_2(A)$.
\begin{prop}~\cite[Theorem 5]{Wh}
Let $A$ be an abelian group with a basis $\{a_i \ | \ i\in I\}$
for a well-ordered index set $I$, and the defining relations
$\{b_{\lambda}\equiv 0\}$. Then the group $\Gamma_2(A)$ is
combinatorial defined by the set of symbolic generators
$\gamma(a_i)$, $i\in I$, and $\gamma(a_i,a_j)$, $i,j\in I$ with
$i<j$ with defining relations $\gamma(b_{\lambda})\equiv0$ and
$\gamma(a_i,b_{\lambda})\equiv0$.
\end{prop}

\begin{example}\label{example2.1}{\rm
We list some examples of the group $\Gamma_2(A)$. The first two examples are direct consequences of the above proposition.
\begin{enumerate}
\item Let $A$ be a free abelian group with a basis $\{a_i \ | \
i\in I\}$ for a well-ordered index set $I$. Then $\Gamma_2(A)$ i s
the free abelian group with a basis given by  $\gamma(a_i)$, $i\in
I$, and $\gamma(a_i,a_j)$, $i,j\in I$ with $i<j$. \item If $A$ is
a cyclic group of finite order $m$ generated by $a_1$, then
$\Gamma_2(A)$ is cyclic of order $m$ or $2m$, according as $m$ is
odd or even, generated by $\gamma(a_1)$. \item Let
$A=\bigoplus_{i\in I} A_i$ for a well-ordered index set $I$.
Then~\cite[Theorem 7]{Wh}
$$
\Gamma_2(A)\cong \bigoplus_{i\in I}\Gamma_2(A_i)\oplus \bigoplus_{\begin{array}{c}
i,j\in I\\
i<j\\
\end{array}}A_i\otimes A_j.
$$
\item For a general abelian group $A$, there is a short exact sequence~\cite[formula (13.8), p.93]{EM}
$$
A\otimes A\rInto^{t} \Gamma_2(A)\rOnto A\otimes\Z/2,
$$
where $t(a\otimes
b)=\gamma(a,b)=\gamma(a+b)-\gamma(a)-\gamma(b)$.\hfill $\Box$
\end{enumerate}

}\end{example}
\subsection{Lower Homology of $K(A,2)$}\label{subsection2.2}
The homology of Eilenberg-MacLane spaces $K(A,n)$ has been studied
in the classical reference~\cite{EM} and other papers. See also
\cite{Breen} for the functorial description of homology groups of
$K(A,2)$ in all dimensions.

\begin{lemma}\label{lemma2.1}~\cite[Theorems 20.5 and 21.1]{EM}
Let $A$ be any abelian group. Then
\begin{enumerate}
\item $H_2(K(A,2))=A$;
\item $H_3(K(A,2))=0$;
\item $H_4(K(A,2))=\Gamma_2(A)$.\hfill $\Box$
\end{enumerate}
\end{lemma}
The homology $H_5(K(A,2))$ becomes a special functor on $A$. Let $R_2(A)=H_5(K(A,2))$. The group $R_2(A)$ for finitely generated abelian group $A$ can be computed as follows~\cite[Section 22]{EM}:
\begin{enumerate}
\item[1)] If $A$ is a cyclic group of order infinite or odd, then $R_2(A)=0$;
\item[2)] If $A=\Z/2^r\Z$ with $r\geq1$, then $R_2(A)\cong \Z/2$.
\item[3)] Let $A=A_1\oplus A_2$. Then $K(A,2)\simeq K(A_1,2)\times K(A_2,2)$. By using K\"unneth theorem together with the fact that $H_1(K(A,2))=H_3(K(A,2))=0$ from Lemma~\ref{lemma2.1}, we have
 $$
 H_5(K(A,2))\cong H_5(K(A_1,2))\oplus H_5(K(A_2,2))\oplus \Tor (H_2(K(A_1,2)),H_2(K(A_2,2))).
 $$
 Thus
 \begin{equation}\label{R_2}
 R_2(A_1\oplus A_2)\cong R_2(A_1)\oplus R_2(A_2)\oplus \Tor(A_1,A_2).
 \end{equation}
\end{enumerate}

Recall the definition of the derived functors in the sense of
Dold-Puppe \cite{DoldPuppe}. Let $F$ be an endofunctor in the
category of abelian groups and $A$ an abelian group. Take a
projective resolution $P_*\to A$. Let $N^{-1}$ be the inverse map
to the normalization map due to Dold-Kan. Then $N^{-1}P_*$ is a
free simplicial resolution of $A$. Then, the $i$-th derived
functor of $F$ applied to the abelian group $A$, is defined as
follows:
$$
L_iF(A)=\pi_i(F(N^{-1}P_*)),\ i\geq 0.
$$
It is a well-known fact that this definition does not depend on a
choice of a projective resolution. In these notations, one has a
natural isomorphism:
$$
R_2(A)=L_1\Gamma_2(A).
$$

\subsection{Whitehead exact sequence.}
Let $X$ be a $(r-1)$-connected $CW$-complex, $r\geq 2.$ There is the following
long exact sequence of abelian groups \cite[Theorem 1]{Wh}:
\begin{equation}\label{white}
\dots \to H_{n+1}(X)\to \Gamma_n(X)\to
\pi_n(X)\buildrel{h_n}\over\to H_n(X)\to \Gamma_{n-1}(X)\to \dots,
\end{equation} where
$\Gamma_n(X)=\mathrm{Im}(\pi_n(\mathrm{sk}_{n-1}(X))\to
\pi_n(\mathrm{sk}_n(X)))$ (here $\mathrm{sk}_i(X)$ is the $i$-th
skeleton of $X$), $h_n$ is the $n$th Hurewicz homomorphism.

The Hurewicz theorem is equivalent to the statement $\Gamma_i(X)=0,\
i\leq r.$ J. H. C. Whitehead computed the term $\Gamma_{r+1}(X)$: In the following theorem, assertion (1) was given in~\cite[Theorem 14]{Wh} and assertion (2) was given the earlier paper~\cite{Wh-1}. According to the remarks in the end of~\cite[Section 14]{Wh}, assertion (2) has been discussed by G. W. Whitehead~\cite{G.W.Wh} as well.

\begin{theorem}\label{theorem2.1}
Let $X$ be a $(r-1)$-connected $CW$-complex with $r\geq 2.$ Then
\begin{enumerate}
\item If $r=2$, then $\Gamma_3(X)\cong \Gamma_2(\pi_2(X))$.
\item If $r>2$, then $\Gamma_{r+1}(X)\cong \pi_r(X)\otimes\Z/2$. \hfill $\Box$
\end{enumerate}
\end{theorem}
The isomorphism $\Gamma_2(\pi_2(X))\to \Gamma_3(X)$ is constructed as follows: Let $\eta\colon S^3\to S^2$ be the Hopf map and let $x\in \pi_2(X)$ be written as the the composite
$$
S^2\rTo^{\tilde x} \mathrm{sk}_2(X)\rInto \mathrm{sk}_3(X).
$$
Then the composite
$$
S^3\rTo^{\eta} S^2\rTo^{\tilde x}\mathrm{sk}_2(X)\rInto \mathrm{sk}_3(X)
$$
defines an element $\eta^*(x)\in \Gamma_3(X)$. According to~\cite[Section 13]{Wh}, the mapping
\begin{equation}\label{eta1}
\eta_1\colon \Gamma_2(\pi_2(X))\to \Gamma_3(X), \ \ \gamma(x)\mapsto \eta^*(x),
\end{equation}
is a well-defined isomorphism of groups. The construction of the isomorphism $\pi_r(X)\otimes\Z/2\to \Gamma_{r+1}(X)$ in assertion (2) is similar.

Recall the description of the functors $\Gamma_{r+2}(X)$ due to
H.-J. Baues \cite{Ba}. Consider the third super-Lie functor
$$
\EuScript L_s^3: \sf Ab\to Ab
$$
defined as
$$
\EuScript L_s^3(A)=im\{A\otimes A\otimes A\buildrel{l}\over\to
A\otimes A\otimes A\}
$$
where
$$
l(a\otimes b\otimes c)=\{a,b,c\}:=a\otimes b\otimes c+b\otimes
a\otimes c-c\otimes a\otimes b-c\otimes b\otimes a,\ a,b,c\in A.
$$
Observe that $\EuScript L_s^3(A)=\ker\{A\otimes
\Lambda^2(A)\buildrel{r}\over\to \Lambda^3(A)\},$ where
$\Lambda^i(A)$ is the $i$th exterior power of $A$ and the map $r$
is given as
$$
r(a\otimes b\wedge c)=a\wedge b\wedge c,\ a,b,c\in A.
$$

Let the complex $X$ be simply connected. Given an abelian group $A$, define  the  map
$$
q: \Gamma_2(A)\otimes A\to \EuScript L_s^3(A)\oplus
\Gamma_2(A)\otimes \mathbb Z/2
$$
by setting
$$
q(\gamma_2(a)\otimes
b)=-\{b,a,a\}+(\gamma_2(a+b)-\gamma_2(a)-\gamma_2(b))\otimes 1,\
a,b\in A.
$$
Define the group $\Gamma_2^2X=\Gamma_2^2(\Gamma_2(\pi_2X)\to
\pi_3X)$ as the pushout:
\begin{equation}\label{push}
\xyma{\Gamma_2(\pi_2(X))\otimes (\pi_2(X)\oplus \mathbb Z/2)
\ar@{->}[r]^{q\oplus id} \ar@{->}[d]^{\eta_1\otimes id} &
\EuScript L_s^3(\pi_2(X))\oplus \Gamma_2(\pi_2(X))\otimes \mathbb Z/2
\ar@{->}[d]\\ \pi_3(X)\otimes ( \pi_2(X) \oplus \mathbb Z/2)
\ar@{->}[r] & \Gamma_2^2(X)}
\end{equation}

\begin{theorem}~\cite[Theorem 3.1]{Ba}\label{theorem2.2}
Let $X$ be a $(r-1)$-connected $CW$-complexes with $r\geq 2$.
\begin{enumerate}
\item[1)] If $r=2$, then there is a natural short exact sequence
$$
0\to \Gamma_2^2(X)\to \Gamma_4(X)\to R_2(\pi_2(X))\to 0.
$$
\item[2)] If $r=3$, then there is a natural exact sequence
$$
0\to \pi_4(X)\otimes \mathbb Z/2\oplus\Lambda^2(\pi_3(X))\to
\Gamma_5(X)\to \Tor(\pi_3(X),\mathbb Z/2)\to 0.
$$
\item[3)] If $r\geq4$,  there is a natural
exact sequence
$$
0\to \pi_{r+1}(X)\otimes \mathbb Z/2\to \Gamma_{r+2}(X)\to
\Tor(\pi_r(X),\mathbb Z/2)\to 0.
$$
\end{enumerate}
\hfill $\Box$
\end{theorem}

Let $A$ be an abelian group. Consider the Hurewicz homomorphism
$$
h_*\colon \pi_*(\Sigma K(A,1))\to \tilde H_*(\Sigma K(A,1))=\tilde H_{*-1}(K(A,1))=\tilde H_{*-1}(A).
$$
Since $H_*(K(A,1))$ is graded commutative ring, the inclusion $$A=H_1(K(A,1))\rInto H_*(K(A,1))$$ induces a ring homomorphism
$$
\lambda\colon \Lambda(A)\longrightarrow H_*(K(A,1)).
$$
By~\cite[Theorem 19.3]{EM}, $\lambda$ is a monomorphism and so we may consider $\Lambda^n(A)\subseteq H_n(K(A,1))=H_n(A)$.
\begin{lemma}
For every abelian group $A$, the Hurewicz image $$\mathrm{Im}(h_{n+1}\colon \pi_{n+1}(\Sigma K(A,1))\to H_n(A))$$ contains the subgroup $\Lambda^n(A)$.
\end{lemma}
\begin{proof}
From the naturality, it suffices to show that the statement holds for a free abelian group $A$.

When $A$ is a free abelian group, then $K(A,1)$ is a (weak) Cartesian product of the circles. Thus $\Sigma K(A,1)$ is a wedge of spheres from the suspension splitting that
$$
\Sigma X\times Y\simeq \Sigma X\vee \Sigma Y\vee \Sigma X\wedge Y
$$
and so the Hurewicz homomorphism induces an epimorphism
$$
h_*\colon \pi_{n+1}(\Sigma K(A,1))\rOnto H_n(A)=\Lambda^n(A)
$$
for a free abelian group $A$. This finishes the proof.
\end{proof}

\subsection{Carlsson construction.} Let $G_*$ be a simplicial group and $X$ a pointed
simplicial set with a base point $*$. Consider the simplicial
group $F^{G_*}(X)$ defined as
$$
F^H(X)_n=\coprod_{x\in X_n}(G_n)_x,
$$
i.e. in each degree $F^{G_*}(X)_n$ is the free product of groups $G_n$
numerated by elements of $X_n$ modulo $(G_n)_*,$ with the
canonical choice of face and degeneracy morphisms. It is proved in
\cite{Carl} that the geometric realization $|F^{G_*}(X)|$ is homotopy
equivalent to the loop space $\Omega(|X|\wedge B|G|)$. The main
example we will consider is the simplicial circle $X=S^1$ with
$$
S_0^1=\{*\},\ S_1^1=\{*,\sigma\},\ S_2^1=\{*, s_0\sigma,
s_1\sigma\},\dots, S_n^1=\{*, x_0,\dots, x_n\},
$$
where $x_i=s_n\dots \hat s_i\dots s_0\sigma$ and the simplicial
group $G_*$ with $G_n=G$ for a given group $G$, with identity
homomorphisms as all face and degeneracy maps. In this case we use
the notation $F^{G}(X)=F^{G_*}(X)$. One has a homotopy equivalence
$$
|F^{G}(S^1)|\simeq \Omega \Sigma K(G,1).
$$ The group $F^{G}(S^1)_n$ is the $n$-fold
free product of $G$:
$$
F^{G}(S^1)_1=G,\ F^{G}(S^1)_2=G*G,\ F^{G}(S^1)_3=G*G*G,\ \dots
$$
We can formally identify $G*G$ with $s_0G*s_1G$, $G*G*G$ with
$s_1s_0G*s_2s_0G*s_2s_1G$, etc, and to define naturally the face
and degeneracy maps:
$$F^{G_*}(S^1):\ \ \ \ldots\ \begin{matrix}\longrightarrow\\[-3.5mm] \ldots\\[-2.5mm]\longrightarrow\\[-3.5mm]
\longleftarrow\\[-3.5mm]\ldots\\[-2.5mm]\longleftarrow \end{matrix}\
G*G*G\ \begin{matrix}\longrightarrow\\[-3.5mm]\longrightarrow\\[-3.5mm]\longrightarrow\\[-3.5mm]\longrightarrow\\[-3.5mm]\longleftarrow\\[-3.5mm]
\longleftarrow\\[-3.5mm]\longleftarrow
\end{matrix}\ G*G\ \begin{matrix}\longrightarrow\\[-3.5mm] \longrightarrow\\[-3.5mm]\longrightarrow\\[-3.5mm]
\longleftarrow\\[-3.5mm]\longleftarrow \end{matrix}\ G.$$

\noindent{\bf Remark.} Consider the second term $F^{G}(S^1)_2=G*G$
and face morphisms $d_0,d_1,d_2: G*G=s_0(G)*s_1(G)\to G$ defined
as
$$
d_0:\begin{cases} s_0(g)\mapsto g\\ s_1(g)\mapsto 1\end{cases},\
d_1:\begin{cases} s_0(g)\mapsto g\\ s_1(g)\mapsto g\end{cases},\
d_2:\begin{cases} s_0(g)\mapsto 1\\ s_1(g)\mapsto g\end{cases}.
$$
There is a natural commutative diagram
\begin{equation}\label{susdi}
\xyma{\pi_3(\Sigma K(G,1)) \ar@{->}[d]^{\simeq} \ar@{^{(}->}[r] &
G\otimes G \ar@{->}[r] \ar@{->}[d]^f & G \ar@{->>}[r] \ar@{->}[d]^{\simeq}& G_{ab}\ar@{->}[d]^{\simeq}\\
\pi_3(\Sigma K(G,1)) \ar@{^{(}->}[r] & (ker(d_1)\cap
ker(d_2))/\mathcal B_2 \ar@{->}[r] & G\ar@{->>}[r] & G_{ab}}
\end{equation}
where $\mathcal B_2$ is the 2-boundary subgroup of $G*G$ and the
map $f$ is defined as\footnote{We use the standard commutator
relations: $[g,h]=g^{-1}h^{-1}gh$}
$$
f(g\otimes h)=[s_0(g)s_1(g)^{-1}, s_0(h)].\mathcal B_2.
$$
There is a natural description of the 2-boundary (see \cite{ElM},
for example):
$$
\mathcal B_2=[ker(d_0),ker(d_1)\cap
ker(d_2)][ker(d_1),ker(d_2)\cap ker(d_0)][ker(d_2),ker(d_0)\cap
ker(d_1)].
$$
Diagram (\ref{susdi}) implies that $f$ is a natural
isomorphism.\\

In the case $G=\mathbb Z,$ the simplicial group $F^{G}(S^1)$ is
identical to the Milnor construction $F(S^1),$ with
$F(S^1)_n$ a free group of rank $n$, for $n\geq 1:$
$$F(S^1):\ \ \ \ldots\ \begin{matrix}\longrightarrow\\[-3.5mm] \ldots\\[-2.5mm]\longrightarrow\\[-3.5mm]
\longleftarrow\\[-3.5mm]\ldots\\[-2.5mm]\longleftarrow \end{matrix}\ F_3\ \begin{matrix}\longrightarrow\\[-3.5mm]\longrightarrow\\[-3.5mm]\longrightarrow\\[-3.5mm]\longrightarrow\\[-3.5mm]\longleftarrow\\[-3.5mm]
\longleftarrow\\[-3.5mm]\longleftarrow
\end{matrix}\ F_2\ \begin{matrix}\longrightarrow\\[-3.5mm] \longrightarrow\\[-3.5mm]\longrightarrow\\[-3.5mm]
\longleftarrow\\[-3.5mm]\longleftarrow \end{matrix}\ \mathbb Z.$$
In this case there is a homotopy equivalence
$$
|F(S^1)|\simeq \Omega S^2
$$
and the construction $F(S^1)$ provides a combinatorial model for
the computation of homotopy groups of the 2-sphere $S^2$. The
construction $F(S^1)$ was studied from the group-theoretical point
of view in \cite{Wu}. It is easy to find the simplicial generators
of the homotopy classes of $\pi_i(F(S^1))=\pi_{i+1}(S^2)$ for
$i=3,4,5$. In order to find these simplicial generators, consider
the sequence of maps between Milnor simplicial constructions
$F(S^4)\to F(S^3)\to F(S^2)\to F(S^1)$ such that the induced
homomorphisms $\mathbb Z=\pi_2(F(S^2))\to \pi_2(F(S^1))=\mathbb Z$
and $\mathbb Z=\pi_3(F(S^3))\to \pi_3(F(S^2))=\mathbb Z/2$ are
epimorphisms and define the homotopy classes of $\pi_3(S^2)$ and
$\pi_4(S^3)$ respectively.
\begin{center}
\begin{tabular}{ccccccccccc}
$F(S^3)_4$ & $\begin{matrix}\longrightarrow\\[-3.5mm]\ldots\\[-2.5mm]\longrightarrow\\[-3.5mm]\longleftarrow\\[-3.5mm]
\ldots\\[-2.5mm]\longleftarrow
\end{matrix}$ & $\mathbb Z$ \\
$\downarrow$& & $\downarrow\eta^2$ \\
$F(S^2)_4$ & $\begin{matrix}\longrightarrow\\[-3.5mm]\ldots\\[-2.5mm]\longrightarrow\\[-3.5mm]\longleftarrow\\[-3.5mm]
\ldots\\[-2.5mm]\longleftarrow
\end{matrix}$ & $F(S^2)_3$ & $\begin{matrix}\longrightarrow\\[-3.5mm]\longrightarrow\\[-3.5mm]\longrightarrow\\[-3.5mm]\longrightarrow\\[-3.5mm]\longleftarrow\\[-3.5mm]
\longleftarrow\\[-3.5mm]\longleftarrow
\end{matrix}$ & $\mathbb Z$ \\
$\downarrow$ & & $\downarrow$ & & $\downarrow\eta$ \\
$F(S^1)_4$ & $\begin{matrix}\longrightarrow\\[-3.5mm]\ldots\\[-2.5mm]\longrightarrow\\[-3.5mm]\longleftarrow\\[-3.5mm]
\ldots\\[-2.5mm]\longleftarrow
\end{matrix}$ & $F(S^1)_3$ & $\begin{matrix}\longrightarrow\\[-3.5mm]\longrightarrow\\[-3.5mm]\longrightarrow\\[-3.5mm]\longrightarrow\\[-3.5mm]\longleftarrow\\[-3.5mm]
\longleftarrow\\[-3.5mm]\longleftarrow
\end{matrix}$ & $F(S^1)_2$ & $\begin{matrix}\longrightarrow\\[-3.5mm] \longrightarrow\\[-3.5mm]\longrightarrow\\[-3.5mm]
\longleftarrow\\[-3.5mm]\longleftarrow \end{matrix}$ & $\mathbb Z$
\end{tabular}
\end{center}
For $n\geq 3,$ the homotopy class of $\pi_n(S^{n-1})$ defined as
$\pi_{n-1}(F(S^{n-2}))$ is generated by
$[s_0(\sigma_{n-2}),s_1(\sigma_{n-2})]$ in $F(S^{n-2})_{n-1}$ (see
\cite{Wu}), where $\sigma_{n-2}$ is a generator of
$F(S^{n-2})_{n-2}=\mathbb Z$. That is, we can define the
simplicial suspension maps $\eta^i: F(S^{i+1})_{i+1}\to
F(S^i)_{i+1}$ by
$$
\eta^i: \sigma_{i+1}\to [s_0(\sigma_i), s_1(\sigma_i)],\ i\geq 1.
$$
Since the generators of $\pi_i(S^2)$ are presented by suspensions
over Hopf fibration for $i=3,4,5$, the simplicial generators of
$\pi_i(F(S^1)),\ i=2,3,4$ are given by the following elements:

\begin{align}
& w_2(x_0,x_1)=[x_0,x_1]\\ &
w_3(x_0,x_1,x_2)=[[x_0,x_2],[x_0,x_1]]\label{susp5el}\\
&
w_4(x_0,x_1,x_2,x_3)=[[[x_0,x_3],[x_0,x_1]],[[x_0,x_2],[x_0,x_1]]].
\end{align}
Here we use the natural notations $x_j:=s_i\dots \hat s_j\dots
s_0(\sigma_1),\ j=0,\dots, i$ for the basis elements in
$F(S^1)_{i+1}$\footnote{One can continue the process of
construction of elements $w_{n+1}(x_0,\dots,x_n)$ by the following
law: $w_{n+1}(x_0,\dots,x_n)=[w_n(x_0,\dots,\hat
x_{n-1},x_n),w_n(x_0,\dots,x_{n-1})].$ In this case, the
16-commutator bracket $w_5(x_0,\dots,x_4)$ corresponds to the
element of order 2 in $\pi_6(S^2)$, but the 32-commutator bracket
$w_6(x_0,\dots, x_5)$ lies in the simplicial boundary subgroup
$\mathcal BF(S^1)_6$ (see \cite{EM}). The construction of a
simplicial generator of the 3-torsion in $\pi_6(S^2)$ is more
tricky: it is possible to find its simplicial representative which
is a product of six brackets of the commutator weight six.}.

\subsection{Spectral sequence} Consider an abelian group $A$ and
its two-step flat resolution
$$
0\to A_1\to A_0\to A\to 0.
$$
By Dold-Kan correspondence, we obtain the following free abelian
simplicial resolution of ~$A$:
$$
N^{-1}(A_1\hookrightarrow A_0):\ \ \ \ \ldots \begin{matrix}\longrightarrow\\[-3.5mm]\longrightarrow\\[-3.5mm]\longrightarrow\\[-3.5mm]\longleftarrow\\[-3.5mm]
\longleftarrow
\end{matrix}\ A_1\oplus s_0(A_0)\ \begin{matrix}\longrightarrow\\[-3.5mm]\longrightarrow\\[-3.5mm]
\longleftarrow \end{matrix}\ A_0.$$

Applying Carlsson construction to the resolution
$N^{-1}(A_1\hookrightarrow A_0)$, we obtain the following
bisimplicial group:
\begin{center}
\begin{tabular}{ccccccccccc}
$F^{N^{-1}(A_1\hookrightarrow A_0)_2}(S^n)_3$ & $\begin{matrix}\longrightarrow\\[-3.5mm]\longrightarrow\\[-3.5mm]\longrightarrow\\[-3.5mm]\longrightarrow\\[-3.5mm]\longleftarrow\\[-3.5mm]
\longleftarrow\\[-3.5mm]\longleftarrow
\end{matrix}$ & $F^{N^{-1}(A_1\hookrightarrow A_0)_2}(S^n)_2$ & $\begin{matrix}\longrightarrow\\[-3.5mm] \longrightarrow\\[-3.5mm]\longrightarrow\\[-3.5mm]
\longleftarrow\\[-3.5mm]\longleftarrow \end{matrix}$ & $N^{-1}(A_1\hookrightarrow A_0)_2$\\
$\downarrow\downarrow\downarrow\uparrow\uparrow$ & & $\downarrow\downarrow\downarrow\uparrow\uparrow$ & & $\downarrow\downarrow\downarrow\uparrow\uparrow$ \\
$F^{A_1\oplus s_0(A_0)}(S^n)_3$ & $\begin{matrix}\longrightarrow\\[-3.5mm]\longrightarrow\\[-3.5mm]\longrightarrow\\[-3.5mm]\longrightarrow\\[-3.5mm]\longleftarrow\\[-3.5mm]
\longleftarrow\\[-3.5mm]\longleftarrow
\end{matrix}$ & $F^{A_1\oplus s_0(A_0)}(S^n)_2$ & $\begin{matrix}\longrightarrow\\[-3.5mm] \longrightarrow\\[-3.5mm]\longrightarrow\\[-3.5mm]
\longleftarrow\\[-3.5mm]\longleftarrow \end{matrix}$ & $A_1\oplus s_0(A_0)$\\
$\downarrow\downarrow\uparrow$ & & $\downarrow\downarrow\uparrow$ & & $\downarrow\downarrow\uparrow$ \\
$F^{A_0}(S^n)_3$ & $\begin{matrix}\longrightarrow\\[-3.5mm]\longrightarrow\\[-3.5mm]\longrightarrow\\[-3.5mm]\longrightarrow\\[-3.5mm]\longleftarrow\\[-3.5mm]
\longleftarrow\\[-3.5mm]\longleftarrow
\end{matrix}$ & $F^{A_0}(S^n)_2$ & $\begin{matrix}\longrightarrow\\[-3.5mm] \longrightarrow\\[-3.5mm]\longrightarrow\\[-3.5mm]
\longleftarrow\\[-3.5mm]\longleftarrow \end{matrix}$ & $A_0$
\end{tabular}
\end{center}
Here the $m$th horizontal simplicial group is Carlsson
construction $F^{N^{-1}(A_1\hookrightarrow A_0)_m}(S^n).$ By the
result of Quillen~\cite{Quillen}, we obtain the following spectral
sequence:
\begin{equation}\label{specseq}
E_{p,q}^2=\pi_q(\pi_p(\Sigma^n K(N^{-1}(A_1\hookrightarrow
A_0),1))\Longrightarrow \pi_{p+q}(\Sigma^n K(A,1)).
\end{equation}

Consider now a non-abelian analog of this spectral sequence, for
$n=1$. Suppose now that a group $G$ is arbitrary, not necessary
abelian. Consider a simplicial resolution of $G$:
$$
G_\bullet\to G,
$$
i.e. $G_\bullet$ is a simplicial group with $\pi_0(G_\bullet)=G,\
\pi_i(G_\bullet)=0,\ i>0.$ Consider the following bisimplicial
group
\begin{center}
\begin{tabular}{ccccccccccc}
$G_2*G_2*G_2$ & $\begin{matrix}\longrightarrow\\[-3.5mm]\longrightarrow\\[-3.5mm]\longrightarrow\\[-3.5mm]\longrightarrow\\[-3.5mm]\longleftarrow\\[-3.5mm]
\longleftarrow\\[-3.5mm]\longleftarrow
\end{matrix}$ & $G_2*G_2$ & $\begin{matrix}\longrightarrow\\[-3.5mm] \longrightarrow\\[-3.5mm]\longrightarrow\\[-3.5mm]
\longleftarrow\\[-3.5mm]\longleftarrow \end{matrix}$ & $G_2$\\
$\downarrow\downarrow\downarrow\uparrow\uparrow$ & & $\downarrow\downarrow\downarrow\uparrow\uparrow$ & & $\downarrow\downarrow\downarrow\uparrow\uparrow$ \\
$G_1*G_1*G_1$ & $\begin{matrix}\longrightarrow\\[-3.5mm]\longrightarrow\\[-3.5mm]\longrightarrow\\[-3.5mm]\longrightarrow\\[-3.5mm]\longleftarrow\\[-3.5mm]
\longleftarrow\\[-3.5mm]\longleftarrow
\end{matrix}$ & $G_1*G_1$ & $\begin{matrix}\longrightarrow\\[-3.5mm] \longrightarrow\\[-3.5mm]\longrightarrow\\[-3.5mm]
\longleftarrow\\[-3.5mm]\longleftarrow \end{matrix}$ & $G_1$\\
$\downarrow\downarrow\uparrow$ & & $\downarrow\downarrow\uparrow$ & & $\downarrow\downarrow\uparrow$ \\
$G_0*G_0*G_0$ & $\begin{matrix}\longrightarrow\\[-3.5mm]\longrightarrow\\[-3.5mm]\longrightarrow\\[-3.5mm]\longrightarrow\\[-3.5mm]\longleftarrow\\[-3.5mm]
\longleftarrow\\[-3.5mm]\longleftarrow
\end{matrix}$ & $G_0*G_0$ & $\begin{matrix}\longrightarrow\\[-3.5mm] \longrightarrow\\[-3.5mm]\longrightarrow\\[-3.5mm]
\longleftarrow\\[-3.5mm]\longleftarrow \end{matrix}$ & $G_0$
\end{tabular}
\end{center}
Again, by the result of Quillen~\cite{Quillen}, we obtain the
following spectral sequence:
\begin{equation}\label{specseq11}
E_{p,q}^2=\pi_q(\pi_p(\Sigma K(G_\bullet,1)))\Longrightarrow
\pi_{p+q}(\Sigma K(G,1)).
\end{equation}
If $G_\bullet$ is a free simplicial resolution, the spectral
sequence (\ref{specseq11}) contains a lot of canonical
differentials of a complicated nature:
$$
\xyma{& & H_4(G)\ar@{->}[ddl] \\ \pi_2(\EuScript L_s^3((G_\bullet)_{ab})\oplus \Gamma_2((G_\bullet)_{ab})\otimes \Z/2)& \pi_2(\Gamma_2((G_\bullet)_{ab})) & H_3(G) \ar@{->}[ddl]\\
\pi_1(\EuScript L_s^3((G_\bullet)_{ab})\oplus \Gamma_2((G_\bullet)_{ab})\otimes \Z/2) & \pi_1(\Gamma_2((G_\bullet)_{ab})) & H_2(G)\\
\EuScript L_s^3(G_{ab})\oplus \Gamma_2(G_{ab})\otimes \Z/2 &
\Gamma_2(G_{ab}) & G_{ab}}
$$

\section{On group $\pi_4(\Sigma K(A,1))$}
\subsection{The group $\pi_4(\Sigma K(A,1))$ for an abelian group $A$}
Let $A$ be an abelian group. Consider the homotopy commutative diagram of fibre sequences
\begin{equation}\label{equation3.1}
\begin{diagram}
\Sigma K(A,1)\wedge K(A,1)&\rTo^{H} &\Sigma K(A,1)&\rTo& K(A,2)=BK(A,1)\\
\dEq&&\dTo&\textrm{pull}&\dTo>{\Delta}\\
\Sigma K(A,1)\wedge K(A,1)&\rTo^{f} &K(A,2)\vee K(A,2)&\rInto& K(A,2)\times K(A,2),\\
\end{diagram}
\end{equation}
where $H$ is the Hopf fibration. Thus we have the following lemma:
\begin{lemma}\label{new1}
There are isomorphisms
$$
\pi_n(\Sigma K(A,1)\wedge K(A,1))\cong \pi_n(\Sigma K(A,1))
\cong \pi_n(K(A,2)\vee K(A,2))
$$
for $n\geq 3$. In particular, $\pi_3(\Sigma K(A,1))\cong \pi_3(\Sigma K(A,1)\wedge K(A,1))\cong A\otimes A$.\hfill $\Box$
\end{lemma}

By Lemma~\ref{lemma2.1}, the lower homology of the wedge $K(A,2)\vee K(A,2)$ are the
following:

\begin{lemma}\label{lemma3.2}
\begin{align*} & H_2(K(A,2)\vee K(A,2))=A\oplus A,\\
& H_3(K(A,2)\vee K(A,2))=0,\\
& H_4(K(A,2)\vee K(A,2))=\Gamma_2(A)\oplus \Gamma_2(A),\\
& H_5(K(A,2)\vee K(A,2))=R_2(A)\oplus R_2(A).\\
\end{align*}\hfill $\Box$
\end{lemma}

\begin{lemma}\label{lemma3.3}
The Hurewicz image
$$
h_n\colon \pi_n(K(A,2)\vee K(A,2))\longrightarrow H_n(K(A,2)\vee K(A,2))
$$
is zero for $n\geq 3$.
\end{lemma}
\begin{proof}
The assertion follows from the commutative diagram
\begin{diagram}
\pi_n(K(A,2)\vee K(A,2))&\rTo^{h_n}& H_n(K(A,2)\vee K(A,2))\\
\dTo&&\dInto\\
\pi_n(K(A,2)\times K(A,2))=0&\rTo^{h_n}&H_n(K(A,2)\times K(A,2)).\\
\end{diagram}
\end{proof}

\begin{theorem}\label{theorem3.1}
Let $A$ be any abelian group. Then there is a natural short exact
sequence
$$
(\Lambda^2(A)\otimes A)^{\oplus 2}\oplus A\otimes A\otimes
\mathbb Z/2\rInto \pi_4(\Sigma K(A,1))\rOnto \Tor(A,A).
$$
Moreover  $(\Lambda^2(A)\otimes A)^{\oplus 2}$ is an (unnatural) summand of $\pi_4(\Sigma K(A,1))$.
\end{theorem}
\begin{proof}
Let $X=K(A,2)\vee K(A,2)$ and let $Y=K(A,2)\times K(A,2)=K(A\oplus A, 2)$. From Lemmas~\ref{lemma3.2} and~\ref{lemma3.3}, there is a short exact sequence
$$
R_2(A)\oplus R_2(A) \rInto \Gamma_4(X)\rOnto \pi_4(X).
$$
The inclusion $j\colon X=K(A,2)\vee K(A,2)\rInto Y=K(A,2)\times K(A,2)$ induces a commutative diagram
\begin{equation}\label{equation3.2}
\begin{diagram}
H_5(X)=R_2(A)\oplus R_2(A)& \rInto& \Gamma_4(X)&\rOnto &\pi_4(X)\\
\dInto>{j_*}&&\dTo>{j_*}&&\dTo>{j_*}\\
H_5(Y)=R_2(A\oplus A)& \rTo^{\cong}_{\phi}& \Gamma_4(Y)&\rTo &\pi_4(Y)=0.\\
\end{diagram}
\end{equation}
By formula~(\ref{R_2}),
$$
H_5(K(A,2)\times K(A,2))=H_5(K(A\oplus A,2))=R_2(A)\oplus R_2(A)\oplus \Tor(A,A)
$$
and so the cokernel of $j_*\colon H_5(X)\to H_5(Y)$ is $\Tor(A,A)$. On the other hand, from Theorem~\ref{theorem2.2}(1), there is a commutative diagram of short exact sequences
\begin{diagram}
\Gamma_2^2(X)&\rInto& \Gamma_4(X)&\rOnto& R_2(\pi_2(X))=R_2(A\oplus A)\\
\dTo&&\dTo>{j_*} &&\dTo>{\cong}\\
\Gamma^2_2(Y)&\rInto& \Gamma_4(Y)&\rOnto^{\psi}& R_2(\pi_2(Y))=R_2(A\oplus A).\\
\end{diagram}
The composite
$$
\psi\circ\phi \colon R_2(A\oplus A)\rTo^{\cong}_{\phi} \Gamma_4(Y)\rOnto_{\psi} R_2(A\oplus A)
$$
is a natural self epimorphism for any abelian group $A$. It is an isomorphism for any finitely generated abelian $A$ and so an isomorphism for any abelian group $A$ by considering the direct limit. Thus $j_*\colon \Gamma_4(X)\to \Gamma_4(Y)$ is an epimorphism and, from diagram~(\ref{equation3.2}), there is a short exact sequence
\begin{equation}\label{equation3.3}
\Gamma_2^2(X)\rInto\pi_4(X)\rOnto \Tor(A,A).
\end{equation}

Let $Z=\Sigma K(A,1)\wedge K(A,1)$ and let $f\colon Z\to X$ be the map in diagram~(\ref{equation3.1}). Consider the commutative diagram
\begin{equation}\label{equation3.4}
\begin{diagram}
\pi_3(Z)\otimes\Z/2=\Gamma_2^2(Z)&\rInto^{\cong}&\Gamma_4(Z)&\rTo&\pi_4(Z)&\rTo& H_4(Z)&\rTo&\Gamma_3(Z)=0\\
\dTo>{f_*}&&\dTo>{f_*}&&f_*\dTo>{\cong}&&&&\\
\Gamma_2^2(X)&\rInto&\Gamma_4(X)&\rOnto&\pi_4(X),&&&&\\
\end{diagram}
\end{equation}
where $\Gamma_2^2(Z)\to \Gamma_4(Z)$ is an isomorphism because its cokernel $R_2(\pi_2(Z))=0$. From the definition~(\ref{push}) of the functor $\Gamma_2^2$,
$$
f_*\colon \Gamma_2^2(Z)\longrightarrow \Gamma_2^2(X)
$$
is a monomorphism with retracting homomorphism $\phi'\colon \Gamma_2^2(X)\to \Gamma_2^2(Z)$. By the short exact sequence~(\ref{equation3.3}), $\Gamma_4(Z)\to \pi_4(Z)$ is a monomorphism and so there is a short exact sequence
\begin{equation}\label{equation3.5}
A\otimes A\otimes\Z/2=\Gamma_4(Z)\rInto \pi_4(Z)\rOnto H_4(Z)=H_3(K(A,1)\wedge K(A,1)).
\end{equation}
Note that $H_1(K(A,1))=A$ and $H_2(K(A,1))=\Lambda^2(A)$. By the K\"unneth theorem, there is a natural short exact sequence
$$
(A\otimes \Lambda^2(A))^{\oplus 2}\rInto H_4(Z)\rOnto^{\psi'} \Tor(A,A).
$$
Consider the composite
$$
\theta_A\colon \Gamma^2_2(X)=\Gamma^2_2(\Gamma_2(A\oplus A)\to A\otimes A)\rInto \pi_4(X)\rTo^{f_*^{-1}}\pi_4(Z)\rOnto H_4(Z)\rOnto^{\psi'}\Tor(A,A),
$$
which is natural on any abelian group $A$. If $A$ is a free abelian group, then $\theta_A=0$. For any abelian group $A$, choose any free abelian group $A_0$ with an epimorphism $g\colon A_0\twoheadrightarrow A$. From the definition~(\ref{push}) of $\Gamma_2^2$,
$$
\Gamma^2_2(g)\colon \Gamma^2_2(\Gamma_2(A_0\oplus A_0)\to A_0\otimes A_0)\longrightarrow \Gamma^2_2(\Gamma_2(A\oplus A)\to A\otimes A)
$$
is an epimorphism. By the naturality of $\theta_A$, we have $\theta_A=0$ because $\theta_{A_0}=0$. Now, from diagram~(\ref{equation3.4}), there is a commutative diagram of natural short exact sequences
\begin{diagram}
A\otimes A\otimes\Z/2&\rEq& A\otimes A\otimes\Z/2&&\\
\dInto>{\circlearrowleft}&&\dInto&&\\
\Gamma_2^2(X)&\rInto&\pi_4(\Sigma K(A,1))&\rOnto&\Tor(A,A)\\
\dOnto&&\dOnto&&\dTo>{\cong}\\
(A\otimes \Lambda^2(A))^{\oplus 2}&\rInto& H_4(Z)&\rOnto&\Tor(A,A).\\
\end{diagram}
It follows that there is a natural (on $A$) isomorphism
$$
\Gamma^2_2(X)\cong A\otimes A\otimes\Z/2 \oplus (A\otimes \Lambda^2(A))^{\oplus 2}.
$$
Since $(A\otimes \Lambda^2(A))^{\oplus 2}$ is an (unnatural) summand of $H_4(Z)$, it is an (unnatural) summand of $\pi_4(\Sigma K(A,1))$. The proof is finished.
\end{proof}

\begin{cor}
Let $p$ be an odd prime integer. Then
$$
\pi_4(\Sigma K(\mathbb Z/p^r,1))=\mathbb Z/p^r
$$
and the Hurewicz homomorphism
$$
\pi_4(\Sigma K(\mathbb Z/p^r,1))\to H_4(\Sigma K(\mathbb Z/p^r,1))
$$
is an isomorphism.
\end{cor}
\begin{proof}
In this case, $A\otimes A\otimes \mathbb{Z}/2=0$. Since $\Z/p^r$ is cyclic, $\Lambda^2(A)\otimes A=0$ and hence the result.
\end{proof}

For completely determining the group $\pi_4(\Sigma K(A,1))$, we have to consider the divisibility problem of the elements in the subgroup $A\otimes A\otimes \Z/2=\Gamma_4(Z)\subseteq \pi_4(\Sigma K(A,1))=\pi_4(Z)$. We solve this problem for any finitely generated abelian group $A$.

\begin{lemma}\label{lemma3.4}
Let $A$ be any abelian group and let $j\colon M(A,1)\to K(A,1)$ be a map such that $j_*\colon H_1(M(A,1))\to H_1(K(A,1))$ is an isomorphism. Then there is an (unnatural) splitting exact sequence
$$
\pi_4(\Sigma M(A,1)\wedge M(A,1))\rInto^{\curvearrowleft}_{(\Sigma j\wedge j)_*} \pi_4(\Sigma K(A,1)\wedge K(A,1))\rOnto^{\curvearrowleft} (A\otimes \Lambda^2(A))^{\oplus 2}.
$$
\end{lemma}
\begin{proof}
Let $X=\Sigma M(A,1)\wedge M(A,1)$ and let $Z=\Sigma K(A,1)\wedge K(A,1)$.
The assertion follows from the commutative diagram of short exact sequences
\begin{diagram}
\Gamma_4(X)&\rInto& \pi_4(X)&\rOnto & H_4(X)=\Tor(A,A)\\
\dTo>{\cong}&&\dTo&&\dInto\\
\Gamma_4(Z)&\rInto& \pi_4(Z)&\rOnto& H_4(Z),\\
\end{diagram}
where the bottom row is short exact by equation~(\ref{equation3.5}).
\end{proof}

Given a finitely generated abelian group $A$, let
$$
A=A_1\oplus\bigoplus_{
\begin{array}{c}
r\geq 1\\
p\textrm{ is a prime}\\
\end{array}}A_{p^r}
$$
be the primary decomposition of $A$, where $A_1$ is torsion free and $A_{p^r}$ is a free $\Z/p^r$-module.

\begin{theorem}\label{theorem3.2}
Let $A$ be any finitely generated abelian group. Let $A=A_2\oplus B$ with $B=A_1\oplus \bigoplus_{p^r\not=2}A_{p^r}$. Then
$$
\begin{array}{rcl}
\pi_4(\Sigma K(A,1))&\cong& \frac{1}{2}(A_2\otimes A_2)\oplus (A_2\otimes B)^{\oplus 2}\oplus B^{\otimes 2}\otimes\Z/2\oplus (A\otimes \Lambda^2(A))^{\oplus 2}\\
 & & \ \ \oplus \Tor(A_2,B)^{\oplus 2}\oplus \Tor(B,B),\\
\end{array}
$$
where $\frac{1}{2}(A_2\otimes A_2)$ is a free $\Z/4$-module with rank of $\dim_{\Z/2}(A_2\otimes A_2)$.
\end{theorem}
\begin{proof}
Let $X=\Sigma M(A,1)\wedge M(A,1)$.
By Lemma~\ref{lemma3.4}, it suffices to show that
$$
\pi_4(X) \cong \frac{1}{2}(A_2\otimes A_2)\oplus (A_2\otimes B)^{\oplus 2}\oplus B^{\otimes 2}\otimes\Z/2\oplus  \Tor(A_2,B)^{\oplus 2}\oplus \Tor(B,B).
$$
Observe that there is a homotopy decomposition
$$
X\simeq \bigvee_{
\begin{array}{c}
r,s\geq0\\
p,q \textrm{ prime}\\
\end{array}}\Sigma M(A_{p^r},1)\wedge M(A_{q^s},1),
$$
where we allow $r,s$ to be $0$ for having the factor $A_1$ to be appeared.
Thus there is a decomposition
\begin{equation}\label{equation3.6}
\pi_4(X)\cong \bigoplus_{
\begin{array}{c}
r,s\geq0\\
p,q \textrm{ prime}\\
\end{array}}\pi_4(\Sigma M(A_{p^r},1)\wedge M(A_{q^s},1)).
\end{equation}
Let $r,s\geq 1$ and let $p$ and $q$ be positive prime integers. From~\cite[Corollary 6.6]{Neisendorfer}, there is a homotopy decomposition
\begin{equation}\label{equation3.7}
\begin{array}{l}
\Sigma M(\Z/p^r,1)\wedge M(\Z/q^s,1)\\
\simeq \left\{
\begin{array}{lcl}
\ast&\textrm{ if }& p\not=q,\\
M(\Z/p^{\min\{r,s\}},3)\vee M(\Z/p^{\min\{r,s\}},4)& \textrm { if }& p=q\textrm{ and }\max\{p^r, q^s\}>2.\\
 \end{array}
 \right.\\
 \end{array}
\end{equation}
By taking $\pi_4$ to above decomposition, we have
\begin{equation}\label{equation3.8}
\pi_4(\Sigma M(\Z/p^r,1)\wedge M(\Z/q^s,1))\cong \Z/p^r\otimes \Z/q^s\otimes\Z/2 \oplus \Tor(\Z/2^r,\Z/2^s)
\end{equation}
if $\max\{p^r,q^s\}>2$. Clearly this formula also holds for the case where $p^r=1$ or $q^s=1$. For the case $p^r=q^s=2$, we claim that
\begin{equation}\label{equation3.9}
\pi_4(\Sigma M(\Z/2,1)\wedge M(\Z/2,1)))=\Z/4.
\end{equation}
Let $Y=\Sigma M(\Z/2,1)\wedge M(\Z/2,1)$. From the short exact sequence
$$
\Gamma_4(Y)=\Z/2\rInto \pi_4(Y)\rOnto H_4(Y)=\Z/2,
$$
the group $\pi_4(Y)=\Z/4$ or $\Z/2\oplus\Z/2$. Suppose that $\pi_4(Y)=\Z/2\oplus\Z/2$. Then there exists an element $\alpha\in \pi_4(Y)$ of order $2$ which has the nontrivial Hurewicz image. Since $\alpha$ is of order $2$, the map $\alpha\colon S^4\to Y$ extends to a map $\tilde\alpha\colon M(\Z/2,4)\to Y$ with
$$
\tilde \alpha_*\colon H_4(M(\Z/2,4))\longrightarrow H_4(Y)
$$
an isomorphism. Let
$$
j\colon M(\Z/2,3)\longrightarrow Y
$$
be the canonical inclusion. Then $j_*\colon H_3(M(\Z/2,3))\to H_3(Y)$ is an isomorphism. Then
$$
(j,\tilde\alpha)\colon  M(\Z/2,3)\vee M(\Z/2,4)\longrightarrow Y
$$
is a homotopy equivalence because it induces an isomorphism on homology, which contradicts that the Steenrod operation $Sq^2\colon H^3(Y;\Z/2)\to H^5(Y;\Z/2)$ is an isomorphism. Thus $\pi_4(Y)=\Z/4$.

Now the assertion follows from decomposition~(\ref{equation3.6}) and the computational formulae~(\ref{equation3.8}) and ~(\ref{equation3.9}).
\end{proof}

\begin{cor}\label{corollary3.2}
Let $A_2$ be any elementary $2$-group. Then there is a natural short exact sequence
$$
(A_2\otimes \Lambda^2(A_2))^{\oplus 2}\rInto \pi_4(\Sigma K(A_2,1))\rOnto \frac{1}{2}(A_2\otimes A_2),
$$
where $\frac{1}{2}(A_2\otimes A_2)$ is a free $\Z/4$-module. Moreover this splits off unnaturally.
\end{cor}
\begin{proof}
By Theorem~\ref{theorem3.1}, there is a natural short exact sequence
$$
(A_2\otimes \Lambda^2(A_2))^{\oplus 2}\rInto \pi_4(\Sigma K(A_2,1))\rOnto \pi_4(\Sigma K(A_2,1))/(A_2\otimes \Lambda^2(A_2))^{\oplus 2}.
$$
By Theorem~\ref{theorem3.2}, the quotient group $\pi_4(\Sigma K(A_2,1))/(A_2\otimes \Lambda^2(A_2))^{\oplus 2}$ is free $\Z/4$-module for any finite dimensional elementary $2$-groups. The assertion follows by taking direct limits.
\end{proof}

\begin{remark}
The summand $\frac{1}{2}(A_2\otimes A_2)$ is sub-quotient functor of $\pi_4(\Sigma K(A,1))$ on $A$ in the following sense. For any abelian group $A$, the $\Z/2$-component $A_2$ is given by the image of
$$
Sq^1_*\colon H_2(A;\Z/2)\longrightarrow H_1(A;\Z/2).
$$
Thus $A\mapsto A_2$ is a sub functor of the identity functor on abelian groups. Then $\pi_4(\Sigma K(A_2,1))$ is a sub functor of $\pi_4(\Sigma K(A,1))$ on $A$ and so
$$
\frac{1}{2}(A_2\otimes A_2)= \pi_4(\Sigma K(A_2,1))/(A_2\otimes \Lambda^2(A_2))^{\oplus 2}
$$
is a sub-quotient functor of $\pi_4(\Sigma K(A,1))$ on $A$.
\end{remark}

\subsection{Applications}
As an application, we compute $\pi_i(M(\Z/p^r,2))$ for $i\leq 4$. By the Hurewicz Theorem, $\pi_2(M(\Z/p^r,2))=\Z/p^r$. From the Whitehead exact sequence~(\ref{white}), we have
\begin{equation}\label{equation3.10}
\Gamma_n(M(\Z/p^r,2))=\pi_n(M(\Z/p^r,2))
\end{equation}
for $r\geq 3$ because $H_i(M(\Z/p^r,2)=0$ for $i\geq3$. It follows directly that
\begin{equation}\label{equation3.11}
\pi_3(M(\Z/p^r,2))=\Gamma_3(M(\Z/p^r,2))=\Gamma_2(\Z/p^r)=\left\{
\begin{array}{lcl}
\Z/p^r&\textrm{ if } & p>2,\\
\Z/2^{r+1}&\textrm{ if }& p=2,\\
\end{array}\right.
\end{equation}
where $\Gamma_2(A)$ is computed in Example~\ref{example2.1}. From Theorem~\ref{theorem2.2} (1), there is a short exact sequence
$$
\Gamma^2_2(M(\Z/p^r,2))\rInto \Gamma_4(M(\Z/p^r,2))\rOnto R_2(\pi_2(M(\Z/p^r,2)))=R_2(\Z/p^r).
$$
From Subsection~\ref{subsection2.2},
$$
R_2(\Z/p^r)=\left\{
\begin{array}{lcl}
0&\textrm{ if } & p>2,\\
\Z/2&\textrm{ if }& p=2.\\
\end{array}\right.
$$
By the definition~(\ref{push}) of the functor $\Gamma_2^2$, the group $\Gamma_2^2(M(\Z/p^r,2))$ is given by the push-out
\begin{diagram}
\Gamma_2(\Z/p^r)\otimes (\Z/p^r\oplus \mathbb Z/2)&\rTo^{q\oplus \mathrm{id}}&\mathcal{L}_s^3(\Z/p^r)\oplus \Gamma_2(\Z/p^r)\otimes\Z/2\\
\cong\dTo>{\eta_1\otimes\mathrm{id}}&&\dTo\\
\pi_3(M(\Z/p^r,2))\otimes(\Z/p^r\oplus \Z/2)&\rTo & \Gamma_2^2(M(\Z/p^r,2)).
\end{diagram}
Thus
$$
\Gamma_2^2(M(\Z/p^r,2))\cong \mathcal{L}_s^3(\Z/p^r)\oplus \Gamma_2(\Z/p^r)\otimes\Z/2.
$$
Since $\Z/p^r$ is cyclic and $\mathcal{L}^3_s(A)$ is isomorphic to the kernel of $A\otimes\Lambda^2(A)\to \Lambda^3(A)$, we have
$$
\mathcal{L}_s^3(\Z/p^r)=0
$$
and so
$$
\Gamma_2^2(M(\Z/p^r,2))=\left\{
\begin{array}{rcl}
0&\textrm{ if }& p>2,\\
\Z/2&\textrm{ if }&p=2.\\
\end{array}\right.
$$
A direct consequence is:
\begin{equation}\label{equation3.12}
\pi_4(\Sigma M(\Z/p^r,2))=0 \textrm{ for } p>2.
\end{equation}
For the case $p=2$, we have the short exact sequence
$$
\Z/2\rInto \pi_4(M(\Z/2^r,2))\rOnto \Z/2.
$$
The remaining problem is to decide whether $\pi_4(M(\Z/2^r,2))$ is equal to $\Z/2\oplus\Z/2$ or $\Z/4$. It has been computed in~\cite{Wu1} that $\pi_4(M(\Z/2,2))=\Z/4$. For $r>1$, the group $\pi_4(M(\Z/2^r,2))$ seems not recorded in references. We are going to determine the group $\pi_4(M(\Z/2^r,2))$ using our methods.

\begin{lemma}\label{lemma3.5}
Let
$$
j\colon M(\Z/2^r,2)\longrightarrow \Sigma K(\Z/2^r,1)
$$
be the canonical map inducing isomorphism on $H_2$. Then
$$
j_*\colon \Gamma_4(M(\Z/2^r,2))\longrightarrow \Gamma_4(\Sigma K(\Z/2^r,1))
$$
is an isomorphism.
\end{lemma}
\begin{proof}
By Theorem~\ref{theorem2.2}(1), there is a commutative diagram of short exact sequences
\begin{diagram}
\Gamma^2_2(M(\Z/2^r,2))=\Z/2&\rInto&\Gamma_4(M(\Z/2^r,2))& \rOnto& R_2(\Z/2^r)=\Z/2\\
\dTo>{j_*}&&\dTo>{j_*}&&\cong \dTo>{j_*} \\
\Gamma^2_2(K(\Z/2^r,2))&\rInto&\Gamma_4(K(\Z/2^r,2))& \rOnto& R_2(\Z/2^r)=\Z/2\\
\end{diagram}
From the Whitehead exact sequence
$$
\Gamma_3(\Z/2^r)=\Z/2^{r+1}\to \pi_3(\Sigma K(\Z/2^r,1))=\Z/2^r\otimes\Z/2^r=\Z/2^r\to 0,
$$
we have
$$
\eta_1\otimes\mathrm{id}\colon \Gamma_2(\Z/2^r)\otimes (\Z/2^r\oplus \mathbb Z/2)\longrightarrow \pi_3(K(\Z/2^r,2))\otimes(\Z/2^r\oplus \Z/2)
$$
is an isomorphism. Similar to the computation of $\Gamma^2_2(M(\Z/2^r,2))$, we have
$$
\Gamma^2_2(K(\Z/2^r,2))=\Z/2
$$
with an isomorphism $j_*\colon \Gamma^2_2(M(\Z/2^r,2))\cong \Gamma^2_2(K(\Z/2^r,2))$. The assertion then follows by $5$-lemma.
\end{proof}

\begin{lemma}\label{lemma3.6}
The group
$$
\Gamma_4(\Sigma K(\Z/2^r,1))=\left\{
\begin{array}{lcl}
\Z/2\oplus\Z/2&\textrm{ if }& r>1,\\
\Z/4&\textrm{ if }&r=1.\\
\end{array}\right.
$$
\end{lemma}
\begin{proof}
Let $Z=\Sigma K(\Z/2^r,1)\wedge K(\Z/2^r,1)$ and let $H\colon Z\to \Sigma K(\Z/2^r,1)$ be the Hopf map. From equation~(\ref{equation3.5}), there is a commutative diagram of exact sequence
\begin{diagram}
\Gamma_4(Z)=\Z/2&\rInto& \pi_4(Z)&\rOnto& H_4(Z)=\Z/2^r\\
\dTo>{H_*}&&\cong\dTo>{H_*}&&\dTo>{H_*}\\
\Gamma_4(\Sigma K(\Z/2^r,1))&\rInto& \pi_4(\Sigma K(\Z/2^r,1))&\rTo& H_4(\Sigma K(\Z/2^r,1)),\\
\end{diagram}
where the bottom row is left exact because $H_5(\Sigma K(\Z/2^r,1))=H_4(\Z/2^r)=0$.

If $r>1$, then $\Gamma_4(Z)$ is a summand of $\pi_4(Z)\cong \pi_4(\Sigma K(\Z/2^r,1))$ by Theorem~\ref{theorem3.2}. Thus $\Gamma_4(Z)=\Z/2$ is also a summand of $\Gamma_4(\Sigma K(\Z/2^r,1))$. It follows that
$$
\Gamma_4(\Sigma K(\Z/2^r,1))=\Z/2\oplus \Z/2 \textrm{ if } r>1.
$$

If $r=1$, by Corollary~\ref{corollary3.2}, $\pi_4(\Sigma K(\Z/2,1))=\Z/4$ and so
$$
\Gamma_4(\Sigma K(\Z/2,1))\cong \pi_4(\Sigma K(\Z/2^r,1))=\Z/4.
$$
The proof is finished.
\end{proof}

Since $\Gamma_4(\Sigma K(\Z/2^r,1))\to \pi_4(\Sigma K(\Z/2^r,1))$ is a monomorphism, from Lemmas~\ref{lemma3.5} and ~\ref{lemma3.6}, we have the following:

\begin{cor}\label{corollary3.3}
Let
$$
j\colon M(\Z/2^r,2)\longrightarrow \Sigma K(\Z/2^r,1)
$$
be the canonical map inducing isomorphism on $H_2$. Then
\begin{enumerate}
\item $\pi_4(M(\Z/2,2))=\Z/4$ and
$$
j_*\colon \pi_4(M(\Z/2,2))\to \pi_4(\Sigma K(\Z/2,1))
$$
is an isomorphism.
\item For $r>1$, $\pi_4(M(\Z/2^r,2))=\Z/2\oplus \Z/2$ and
$$
j_*\colon \pi_4(M(\Z/2^r,2))=\Z/2\oplus \Z/2\to \pi_4(\Sigma K(\Z/2^r,1))=\Z/2\oplus \Z/2^r
$$
is a monomorphism. \hfill $\Box$
\end{enumerate}
\end{cor}

Note that $M(\Z/2,2)=\Sigma \RP^2$ and $\Sigma K(\Z/2,1)=\Sigma \RP^\infty$ with the canonical inclusion $j\colon \Sigma \RP^2\hookrightarrow \Sigma \RP^\infty$. A consequence of Corollary~\ref{corollary3.3} (1) on the suspended projective spaces are as follows.

\begin{cor}\label{corollary3.4}
Let $j\colon \Sigma \RP^2\to \Sigma \RP^n$ be the canonical inclusion with $3\leq n\leq \infty$.
\begin{enumerate}
\item For $4\leq n\leq\infty$, $j_*\colon \pi_4(\Sigma\RP^2)=\Z/4\to \pi_4(\Sigma \RP^n)$ is an isomorphism.
\item For $n=3$, $j_*\colon \pi_4(\Sigma \RP^2)=\Z/4\to \pi_4(\Sigma \RP^3)$ is a splitting monomorphism. Moreover
$$
\pi_4(\Sigma \RP^3)\cong \pi_4(\Sigma\RP^2)\oplus\Z=\Z/4\oplus\Z.
$$
\end{enumerate}
\end{cor}
\begin{proof}
Assertion (1) and the first part of assertion (2) are direct consequences of Corollary~\ref{corollary3.3}. For the second part of assertion (2), notice that $\RP^3=SO(3)$. From the commutative diagram
\begin{diagram}
\pi_4(\Sigma SO(3)\wedge SO(3))&\rTo^{\cong}&\pi_4(\Sigma \RP^\infty\wedge \RP^\infty)\\
\dTo>{H_*}&&\cong\dTo>{H_*}\\
\pi_4(\Sigma SO(3))&\rTo&\pi_4(\Sigma \RP^\infty),\\
\end{diagram}
we have
$$
\begin{array}{rcl}
\pi_4(\Sigma SO(3))&\cong&\pi_4(\Sigma SO(3)\wedge SO(3))\oplus \pi_4(BSO(3))\\
&\cong& \pi_4(\Sigma \RP^\infty)\oplus \pi_3(SO(3))\\
&\cong&\Z/4\oplus\Z\\
\end{array}
$$
and hence the result.
\end{proof}

Another consequence is as follows:

\begin{cor}\label{corollary3.5}
Let $\Sigma_3$ be the third symmetric group. Then
$\pi_4(\Sigma K(\Sigma_3,1))=\mathbb Z/12.$
\end{cor}
\begin{proof}
Recall that the integral homology groups of $\Sigma_3$ are 4-periodic with the
following initial terms:
$$
H_1(\Sigma_3)=\mathbb Z/2,\ H_2(\Sigma_3)=0,\ H_3(\Sigma_3)=\mathbb
Z/6,\ H_4(\Sigma_3)=0.
$$
Let $X=\Sigma K(\Sigma_3,1)$. The Whitehead exact sequence has the following form:
$$
\Gamma_4(X)\rInto\pi_4(X)\rTo H_3(\Sigma_3)=\Z/6\rTo \Gamma_3(X)=\Z/4\rOnto \pi_3(X)=\Z/2.
$$
The inclusion $\Sigma_2=\Z/2\to \Sigma_3$ induces an isomorphism
$$
\pi_i(\Sigma K(\Z/2,1))=\Z/2\rTo^{\cong} \pi_i(\Sigma K(\Sigma_3,1))=\Z/2
$$
for $i=2,3$.  By Theorem~\ref{theorem2.2} (1) together with Lemma~\ref{lemma3.6}, the inclusion $\Sigma_2=\Z/2\to\Sigma_3$ induces an
$$
\Gamma_4(\Sigma K(\Z/2,1))=\Z/4\rTo^{\cong}\Gamma_4(\Sigma K(\Sigma_3,1))
$$
and hence the result.
\end{proof}

\section{On group $\pi_5(\Sigma K(A,1))$}
\subsection{Some Properties of the Functor $A\mapsto \pi_5(\Sigma K(A,1))$}
From Hopf fibration
$$
\Sigma K(A,1)\wedge K(A,1)\rTo \Sigma K(A,1)\rTo K(A,2),
$$
it suffices to compute $\pi_5(\Sigma K(A,1)\wedge K(A,1))$. Let $Z=\Sigma K(A,1)\wedge K(A,1)$. Since $Z$ is $2$-connected, from Theorem~\ref{theorem2.2}(2),  there are
natural exact sequences
\begin{equation}\label{e433}
\xyma{& \pi_4(\Sigma K(A,1))\otimes \Z/2\oplus \Lambda^2(A\otimes
A) \ar@{^{(}->}[d]\\ H_6(Z)\ar@{->}[r] & \Gamma_5(Z)\ar@{->}[r]
\ar@{->>}[d] & \pi_5(Z)\ar@{->>}[r] & H_5(Z),\\
& \Tor(A\otimes A,\mathbb Z/2)}
\end{equation}
where $\pi_5(Z)\to H_5(Z)$ is onto by equation~(\ref{equation3.5}). The group $\pi_4(\Sigma K(A,1))$ has been determined by Theorems~\ref{theorem3.1} and~\ref{theorem3.2}.
\begin{prop}\label{proposition4.1}
Let $A$ be a free abelian group. Then there is a natural short exact sequence
$$
(\Lambda^2(A)\otimes A\otimes \Z/2)^{\oplus 2}\oplus A^{\otimes 2}\otimes \Z/2\oplus \Lambda^2(A\otimes A)\hookrightarrow\pi_5(\Sigma
K(A,1))
\twoheadrightarrow  (\Lambda^3(A)\otimes A)^{\oplus 2}\oplus
\Lambda^2(A)^{\otimes 2}.
$$
\end{prop}
\begin{proof}
Since $A$ is a free abelian group, the Hurewicz homomorphism $h_*\colon \pi_*(Z)\to \tilde H_*(Z)$ is onto because $Z$ is a wedge of spheres. Thus there is a short exact sequence
$$
\Gamma_5(Z)\rInto \pi_5(Z)\rOnto H_5(Z).
$$
By Theorem~\ref{theorem3.1}, $\pi_4(Z)\cong (\Lambda^2(A)\otimes A)^{\oplus 2}\oplus A\otimes A\otimes\Z/2$ for a free abelian group $A$. The assertion follows from diagram~(\ref{e433}).
\end{proof}

\begin{prop}\label{proposition4.2}
If $A$ is a torsion abelian group with the property that $2\colon A\to A$ is an isomorphism, then there is a natural short exact sequence
$$
\Lambda^2(A\otimes A)\rInto \pi_5(\Sigma K(A,1))\rOnto H_4(K(A,1)\wedge K(A,1)).
$$
\end{prop}
\begin{proof}
It suffices to show that the Hurewicz homomorphism $h_*\colon \pi_6(Z)\to H_6(Z)$ is onto. We may assume that $A$ is finitely generated because we can take direct limit for general case whence the finitely generated case is proved. Then $A$ is a direct sum of the primary $p$-torsion groups $\Z/p^r$ for some $r\geq 1$ and odd primes $p$. According to~\cite{Harris-Kuhn}, there is homotopy decomposition
$$
\Sigma K(\Z/p^r,1)\simeq X_1\vee\cdots\vee X_{p-1},
$$
where $\bar H_q(X_i;\Z)\not=0$ if and only if $q\equiv 2i\mod{2p-2}$. Together with the decomposition formula~(\ref{equation3.7}) for the smash product of Moore spaces, up to $6$-skeleton, $\Sigma K(A,1)\wedge K(A,1)$ is homotopy equivalent to a wedge of spheres and Moore spaces. It follows that the Hurewicz homomorphism
$$
\pi_6(\Sigma K(A,1)\wedge K(A,1))\longrightarrow H_6(\Sigma K(A,1)\wedge K(A,1))
$$
is onto and hence the result.
\end{proof}

From the above proof, we also have the following:

\begin{prop}\label{proposition4.3}
Let $A$ be any abelian group. Let $\Z_{\frac{1}{2}}=\{\frac{m}{2^r}\in \mathbb{Q}\ | \ m\in\Z,\ r\geq0\}$. Then there is natural short exact sequence
$$
\Lambda^2(A\otimes A)\otimes \Z_{\frac{1}{2}} \rInto \pi_5(\Sigma K(A,1))\otimes \Z_{\frac{1}{2}}\rOnto H_4(K(A,1)\wedge K(A,1))\otimes \Z_{\frac{1}{2}}.
$$
\hfill $\Box$
\end{prop}

For computing the group $\pi_5(\Sigma K(A,1))$, as one see from the above, the tricky part is the $2$-torsion. Whence $A$ contains $2$-torsion summands, the Hurewicz homomorphism $\pi_6(Z)\to H_6(Z)$ is no longer epimorphism in general and so $\Gamma_5(Z)\to \pi_5(Z)$ is not a monomorphism in general. Also the group $\pi_5(Z)$ in diagram~(\ref{e433}) admits non-trivial extension. The computation of the group $\pi_5(\Sigma K(A,1))$ for finitely generated abelian groups $A$ can be given by the following steps:
\begin{enumerate}
\item[\textbf{Step 1.}] Take a primary decomposition of $A$ and write $K(A,1)$ as a product of copies of $S^1=K(\Z,1)$ and $K(\Z/p^r,1)$.
\item[\textbf{Step 2.}] By using the fact that $\Sigma X\times Y\simeq \Sigma X\vee \Sigma Y\simeq \Sigma X\wedge Y$ for any spaces $X$ and $Y$, one gets
$$
\begin{array}{rcl}
\Sigma(X_1\times X_2)\wedge (X_1\times X_2)&\simeq& \Sigma (X_1\vee X_2\vee X_1\wedge X_2)\wedge (X_1\vee X_2\vee X_1\wedge X_2)\\
&\simeq& \Sigma (X_1^{\wedge 2}\vee X_2^{\wedge 2}\vee X_1^{\wedge 2}\wedge X_2^{\wedge 2}\vee \\
& &\ \ \bigvee\limits^2X_1\wedge X_2\vee \bigvee\limits^2X_1^{\wedge 2}\wedge X_2\vee\bigvee\limits^2X_1\wedge X_2^{\wedge 2}).\\
\end{array}
$$
From this, $\Sigma K(A,1)\wedge K(A,1)$ is then homotopy equivalent to a wedge of the spaces in the form
$$
X=\Sigma^m K(\Z/p_1^{r_1},1)\wedge K(\Z/p_2^{r_2},1)\wedge \cdots \wedge  K(\Z/p_t^{r_t},1)
$$
with $m+t\geq 3$ and $m\geq1$.
\item[\textbf{Step 3.}] By applying the Hilton-Milnor Theorem, we have
 $$
 \begin{array}{rcl}
 \Omega(\Sigma X\vee \Sigma Y)&\simeq& \Omega \Sigma X \times \Omega \Sigma Y \times \Omega \Sigma ((\Omega\Sigma X)\wedge(\Omega\Sigma Y))\\
 &\simeq&\Omega\Sigma X\times \Omega\Sigma Y\times \Omega\Sigma\left(\bigvee_{i,j=1}^\infty X^{\wedge i}\wedge Y^{\wedge j}\right).\\
 \end{array}
$$
Thus
$$
\pi_n(\Sigma X\vee \Sigma Y)\cong \pi_n(\Sigma X)\oplus \pi_n(\Sigma Y)\oplus \pi_n\left(\bigvee_{i,j=1}^\infty \Sigma X^{\wedge i}\wedge Y^{\wedge j}\right).
$$
Note that the connectivity of $X^{\wedge i}\wedge Y^{\wedge j}$ tends to $\infty$ as $i,j\to \infty$. By repeating the above procedure, $\pi_n(\Sigma K(A,1)\wedge K(A,1))$ is isomorphism to a direct sum of the groups
$\pi_n(X)$
with $X$ given in the form above.
\end{enumerate}
Notice that
$$
\Sigma^m K(\Z/p_1^{r_1},1)\wedge K(\Z/p_2^{r_2},1)\wedge \cdots\wedge  K(\Z/p_t^{r_t},1)\simeq\ast
$$
if the primes $p_i\not=p_j$ for some $i\not=j$. Thus we only need to compute
$$
\pi_5(\Sigma^m K(\Z/p^{r_1},1)\wedge K(\Z/p^{r_2},1)\wedge \cdots\wedge  K(\Z/p^{r_t},1))
$$
for a prime $p$. If $t=0$, the homotopy group $\pi_5(S^m)$ is known by $\pi_5(S^3)=\pi_5(S^4)=\Z/2$ and $\pi_5(S^5)=\Z$. For an odd prime $p$, this homotopy group can be determined by Proposition~\ref{proposition4.2}. The rest work in this section is of course to compute
$
\pi_5(\Sigma^m K(\Z/2^{r_1},1)\wedge K(\Z/2^{r_2},1)\wedge \cdots\wedge  K(\Z/2^{r_t},1))
$
with $m+t\geq 3$. When $m+t\geq 5$, we have
$$
\pi_5(X)=\left\{
\begin{array}{lcl}
0&\textrm{ if }& m+t>5\\
\Z/2^{\min\{r_1,\ldots,r_t\}}&\textrm{ if }& m+t=5 \textrm{ with } t\geq1\\
\end{array}\right.
$$
for $X=\Sigma^m K(\Z/2^{r_1},1)\wedge K(\Z/2^{r_2},1)\wedge \cdots\wedge  K(\Z/2^{r_t},1)$. The first less obvious case is $m+t=4$, which will be discussed in the next subsection.

\subsection{The Group $\pi_5(\Sigma^m K(\Z/2^{r_1},1)\wedge K(\Z/2^{r_2},1)\wedge \cdots\wedge  K(\Z/2^{r_t},1))$ for $m+t=4$ and $m,t\geq1$.}
We first consider the case $t=1$.
\begin{lemma}\label{lemma4.1}
The Hurewicz homomorphism
$$
h_5\colon \pi_5(\Sigma^2 K(\Z/2^r,1))\to H_5(\Sigma^2 K(\Z/2^r,1))
$$
is onto for any $r\geq 1$.
\end{lemma}
\begin{proof}
Let $X=\Sigma^2 K(\Z/2^r,1)$. Consider the Whitehead exact sequence
$$
\pi_5(X)\rTo^{h_5} H_5(X)\to \Gamma_4(X)=\Z/2\to \pi_4(X)\to H_4(X)=0.
$$
Thus the Hurewicz homomorphism $h_5$ is onto if and only if  $\pi_4(X)\not=0$.

Let $f\colon S^3\to X$ be a map representing the generator for $\pi_3(X)=\Z/2^r$. From the remark to Theorem~\ref{theorem2.1}, $\pi_4(X)=0$ if and only if the composite
$$
S^4\rTo^{\eta} S^3\rTo^{f} X
$$
is null homotopic, if and only if the map $f\colon S^3\to X$ extends to a map $\tilde f\colon \Sigma\CP^2\to X$ because $\Sigma\CP^2$ is the homotopy cofibre of $\eta\colon S^4\to S^3$.

Suppose that there exists a map $\tilde f\colon\Sigma\CP^2\to X$ such that $\tilde f|_{S^4}=f$. By taking mod $2$ cohomology, there is commutative diagram
\begin{diagram}
H^5(\Sigma\CP^2;\Z/2)&\lTo^{\tilde f^*}&H^5(X;\Z/2)=\Z/2\\
\cong\uTo>{Sq^2}& &\uTo>{Sq^2}\\
H^3(\Sigma\CP^2;\Z/2)=H^3(S^3;\Z/2)&\lTo^{\tilde f^*=f^*}_{\cong}&H^3(X;\Z/2)=\Z/2.\\
\end{diagram}
It follows that
$$
Sq^2\colon H^3(X;\Z/2)\longrightarrow H^5(X;\Z/2)
$$
is an isomorphism. On the other hand, from the fact that $X=\Sigma^2K(\Z/2^r,1)$ and $Sq^2\colon H^1(K(\Z/2^r,1);\Z/2)\to H^3(K(\Z/2^r,1))$ is zero,  $Sq^2\colon H^3(X;\Z/2)\to H^5(X;\Z/2)$ is zero. This gives a contradiction. The assertion follows.
\end{proof}

\begin{prop}\label{proposition4.4}
$\pi_5(\Sigma^3 K(\Z/2^r,1))=\Z/2$ for $r\geq1$.
\end{prop}
\begin{proof}
Let $X=\Sigma^3 K(\Z/2^r,1)$. Consider the Whitehead exact sequence
$$
\pi_6(X)\rTo^{h_6} H_6(X)\to \Gamma_5(X)=\Z/2\to \pi_5(X)\to H_5(X)=0.
$$
By Lemma~\ref{lemma4.1}, $h_6\colon \pi_6(X)\to H_6(X)$ is onto. Thus $\pi_5(X)\cong \Gamma_5(X)=\Z/2$.
\end{proof}

Now we consider the case $t=2$.
\begin{prop}\label{proposition4.5}
Let $r_1,r_2\geq1$. Then
$$
\pi_5(\Sigma^2 K(\Z/2^{r_1},1)\wedge K(\Z/2^{r_2},1))=\left\{
\begin{array}{lcl}
\Z/2\oplus \Z/2^{\min\{r_1,r_2\}}& \textrm{ if } & \max\{r_1,r_2\}>1,\\
\Z/4&\textrm{ if }& r_1=r_2=1.\\
\end{array}
\right.
$$
\end{prop}
\begin{proof}
Let $X=\Sigma^2 K(\Z/2^{r_1},1)\wedge K(\Z/2^{r_2},1)$. By Lemma~\ref{lemma4.1}, there exists a map
$$
f_i\colon S^5\longrightarrow \Sigma^2 K(\Z/2^{r_i}, 1), \ \ i=1,2,
$$
which induces an epimorphism
$$
f_{i\ast}\colon H_5(S^5)\rOnto H_5(\Sigma^2 K(\Z/2^{r_i},1).
$$
Let  $j\colon Y=\Sigma^2 M(\Z/2^{r_1},1)\wedge M(\Z/2^{r_2},1)\hookrightarrow X$ be the canonical inclusion. Then the map
$$
f\colon Y\vee S^5\wedge K(\Z/2^{r_2},1)\vee K(\Z/2^{r_1},1)\wedge S^5\rTo^{(j, f_1\wedge \mathrm{id},\mathrm{id}\wedge f_2)} X
$$
induces an isomorphism on $H_j(\ ;\Z/2)$ for $j\leq 6$. Thus
$$
f_*\colon \pi_k\left( Y\vee S^5\wedge K(\Z/2^{r_2},1)\vee K(\Z/2^{r_1},1)\wedge S^5\right)\longrightarrow \pi_k(X)
$$
is an isomorphism for $k\leq 5$. Note that
$$
\pi_k(S^5\wedge K(\Z/2^{r_2},1))=\pi_k(K(\Z/2^{r_1},1)\wedge S^5)=0
$$
for $k\leq 5$. Thus
$$
j_*\colon \pi_k(Y)\to \pi_k(X)
$$
is an isomorphism for $k\leq 5$. In particular, $\pi_5(Y)\cong \pi_5(X)$.

If $\max{r_1,r_2}>1$, from decomposition~(\ref{equation3.7}), we have
$$
\Sigma^2 M(\Z/2^{r_1},1)\wedge M(\Z/2^{r_2},1)\simeq M(\Z/2^{\min\{r_1,r_2\}},4)\vee M(\Z/2^{\min\{r_1,r_2\}},5)
$$
and so
$$
\pi_5(Y)\cong \pi_5(M(\Z/2^{\min\{r_1,r_2\}},4))\oplus \pi_5(M(\Z/2^{\min\{r_1,r_2\}},5))=\Z/2\oplus \Z/2^{\min\{r_1,r_2\}}.
$$

Consider the case $r_1=r_2=1$. From formula~(\ref{equation3.9}) and the Freudenthal Suspension Theorem,
$$
\pi_5(\Sigma^2 M(\Z/2,1)\wedge M(\Z/2,1))\cong \pi_4(\Sigma M(\Z/2,1)\wedge M(\Z/2,1))\cong \Z/4.
$$
The proof is finished.
\end{proof}
The last case is $t=3$.
\begin{prop}\label{proposition4.6}
Let $r_1,r_2,r_3\geq 1$ and let $r=\min\{r_1,r_2,r_3\}$. Then
$$
\pi_5(\Sigma K(\Z/2^{r_1},1)\wedge K(\Z/2^{r_2},1)\wedge K(\Z/2^{r_3},1))=\left\{
\begin{array}{lcl}
\Z/2\oplus \Z/2^r&\textrm{ if }& \max\{r_1,r_2,r_3\}>1,\\
\Z/2\oplus\Z/2&\textrm{ if }& r_1=r_2=r_3=1.\\
\end{array}\right.
$$
\end{prop}
\begin{proof}
Let $X=\Sigma K(\Z/2^{r_1},1)\wedge K(\Z/2^{r_2},1)\wedge K(\Z/2^{r_3},1)$. Let $f_1$ be the composite
$$
S^6 \rTo^{g} \Sigma^3 K(\Z/2^{r_1},1)\cong \Sigma K(\Z/2^{r_1},1)\wedge S^1\wedge S^1\rInto X,
$$
where $g$ is a map which induces epimorphism on $H_6(\ )$ by Lemma~\ref{lemma4.1}. Similarly, we have the maps
$$
f_i\colon S^6\longrightarrow \Sigma K(\Z/2^{r_1},1)\wedge K(\Z/2^{r_2},1)\wedge K(\Z/2^{r_3},1),\ \ i=2,3,
$$
by replacing $K(\Z/2^{r_1},1)$ by $K(\Z/2^{r_i},1)$. Let $Y=\Sigma M(\Z/2^{r_1},1)\wedge M(\Z/2^{r_2},1)\wedge M(\Z/2^{r_3},1)$ and let
$
j\colon Y\hookrightarrow X
$
be the canonical inclusion. The map
$$
Y\vee S^6\vee S^6\vee S^6\rTo^{(j,f_1,f_2,f_3)} X
$$
induces an isomorphism on $H_k(\ ;\Z/2)$ for $k\leq 6$ and so
$$
(j,f_1,f_2,f_3)_*\colon \pi_5(Y\vee S^6\vee S^6\vee S^6)=\pi_5(Y)\longrightarrow \pi_5(X)
$$
is an isomorphism.

If $\max\{r_1,r_2,r_3\}>1$, from decomposition~(\ref{equation3.7}),
$$
Y\simeq M(\Z/2^r,4)\vee M(\Z/2^r,5)\vee M(\Z/2^r,5)\vee M(\Z/2^r,6)
$$
and so
$$
\pi_5(Y)=\Z/2\oplus\Z/2^r\oplus\Z/2^r.
$$

If $r_1=r_2=r_3=1$, there is a homotopy decomposition~\cite[Corollary 3.7]{Wu1}
$$
\Sigma \RP^2\wedge\RP^2\wedge \RP^2\simeq \Sigma \CP^2\wedge \RP^2\vee \Sigma^4\RP^2\vee \Sigma^4\RP^2.
$$
By~\cite[Lemma 6.34 (2)]{Wu1},
$$
\pi_5(\Sigma \CP^2\wedge \RP^2)=0
$$
and so
$$
\pi_5(\Sigma \RP^2\wedge\RP^2\wedge \RP^2)=\Z/2\oplus\Z/2,
$$
which finishes the proof.
\end{proof}

\begin{remark}{\rm
For the case $X=\Sigma K(\Z/2,1)\wedge K(\Z/2,1)\wedge K(\Z/2,1)$, the Hurewicz homomorphism
$$
\pi_5(X)=\Z/2\oplus\Z/2\longrightarrow H_5(X)=\Z/2\oplus\Z/2
$$
is an isomorphism and so, in the Whitehead exact sequence,
$$
H_6(X)\longrightarrow \Gamma_5(X)=\Z/2
$$
is onto. This gives an example that the morphism $H_6(Z)\to \Gamma_5(Z)$ in diagram~(\ref{e433}) may not be zero, which is the only example in the case $m+t=4$. More examples will be shown up in the case $m+t=3$ in the next subsections.}
\end{remark}

\subsection{The Group $\pi_5(\Sigma K(\Z/2^r,1))\cong\pi_5(\Sigma K(\mathbb Z/2^r,1)\wedge K(\Z/2^r,1))$.}

\begin{lemma}\label{lemma4.2}
Let $X=\Sigma K(\Z/2^{r},1)\wedge K(\Z/2^{r},1)$ with $r\geq 1$.
Then mod $2$ Hurewicz homomorphism
$$
\pi_6(X)\rTo^{h_6} H_6(X)\rTo H_6(X;\Z/2)
$$
is zero.
\end{lemma}

\begin{proof}
Recall that the mod $2$ cohomology ring
$$
H^*(K(\Z/2^r,1);\Z/2)\cong E(u_1)\otimes P(u_2)
$$
with the $r\,$th Bockstein $\beta^r(u_1)=u_2$. Let $x_i$ (and $y_i$) denote the basis for $H_i(K(\Z/2^{r};\Z/2)$. The Steenrod operations and the Bockstein on lower homology are given by
$$
\begin{array}{rclcrcl}
Sq^2_*x_4&=&x_2&\quad& Sq^2_*y_4&=&y_2\\
\beta_{r}(x_4)&=&x_3&\quad& \beta_{r}y_4&=&y_3\\
\beta_{r}(x_2)&=&x_1&\quad&\beta_{r}y_2&=&y_1.\\
\end{array}
$$
The $\Z/2$-vector space $s^{-1}\tilde H_k(X;\Z/2)$ with $k\leq 6$ has a
basis given by the table
$$
\left(
\begin{array}{lccccc}
k=6&&x_1y_4& x_2y_3& x_3y_2& x_4y_1\\
5&&x_1y_3&x_2y_2&x_3y_1&\\
4&&      &     &x_1y_2&x_2y_1\\
3&&      &     &      &x_1y_1\\
\end{array}
\right)
$$
Let
$\alpha\in H_6(X;\Z/2)$ be a spherical class. Then
$$
s^{-1}\alpha=\epsilon_1x_1y_4+\epsilon_2x_2y_3+\epsilon_3x_3y_2+\epsilon_4x_4y_1
$$
for some $\epsilon_i\in\Z/2$. Observe that for any spherical class,
$$\beta_s(\alpha)=Sq^t_*(\alpha)=0$$ for any $s,t\geq 1$.  By applying $Sq^2_*$ to $\alpha$, we have
$$
\begin{array}{rcl}
0&=&Sq^2_*(s^{-1}\alpha)\\
&=&\epsilon_1Sq^2_*(x_1y_4)+\epsilon_2Sq^2_*(x_2y_3)+\epsilon_3Sq^2_*(x_3y_2)+\epsilon_4Sq^2_*(x_4y_1)\\
&=&\epsilon_1x_1y_2+0+0+\epsilon_4x_2y_1\\
\end{array}
$$
in $s^{-1}H_4(X;\Z/2)$. Thus
\begin{equation}\label{equation-lemma4.2.1}
\epsilon_1=\epsilon_4=0
\end{equation}
By applying the Bockstein $\beta_{r}$ to $\alpha$, we have
$$
\begin{array}{rcl}
0&=&\beta_{r}(s^{-1}\alpha)\\
&=&\epsilon_1\beta_r(x_1y_4)+\epsilon_2\beta_r(x_2y_3)+\epsilon_3\beta_r(x_3y_2)+\epsilon_4\beta_r(x_4y_1)\\
&=&\epsilon_1x_1y_3+\epsilon_2x_1y_3+\epsilon_3x_3y_1+\epsilon_4x_3y_1\\
&=&(\epsilon_1+\epsilon_2)x_1y_3+(\epsilon_3+\epsilon_4)x_3y_1\\
\end{array}
$$
and so
$$
\epsilon_1+\epsilon_2=\epsilon_3+\epsilon_4=0.
$$
Together with equation~(\ref{equation-lemma4.2.1}), we have $\epsilon_i=0$ for $1\leq i\leq 4$. Thus $\alpha=0$ and hence the result.
\end{proof}

\begin{theorem}\label{theorem4.1}
$\pi_5(\Sigma K(\Z/2,1))\cong \pi_5(\Sigma K(\Z/2,1)\wedge K(\Z/2,1))=\Z/2\oplus\Z/2$.
\end{theorem}
\begin{proof}
Let $X=\Sigma K(\Z/2,1)\wedge K(\Z/2,1)$. Notice that
$$
H_6(X)=\Z/2\oplus\Z/2\cong H_6(X;\Z/2).
$$
From diagram~(\ref{e433}), there is an exact sequence
$$
H_6(X)=\Z/2\oplus\Z/2\rInto \Gamma_5(X)\rTo \pi_4(X)\rOnto H_5(X)=\Z/2\oplus\Z/2.
$$
By Corollary~\ref{corollary3.2},
$$
\pi_4(X)\cong \pi_4(\Sigma K(\Z/2,1))\cong\Z/4.
$$
From Theorem~\ref{theorem2.2}(2), there is a short exact sequence
$$
\pi_4(X)\otimes \Z/2\oplus \Lambda^2(\pi_3(X))=\Z/2\rInto \Gamma_5(X)\rOnto \Tor(\pi_3(X),\Z/2)=\Z/2.
$$
Thus the group $\Gamma_5(X)$ is of order $4$. It follows that the monomorphism
$$
H_6(X)=\Z/2\oplus\Z/2 \rInto \Gamma_4(X)
$$
is an isomorphism and hence the result.
\end{proof}

\begin{lemma}\label{lemma4.3}
Let $r_1,r_2\geq 1$ with $\max\{r_1,r_2\}>1$. Then there is a short exact sequence
$$
\Z/2\oplus\Z/2\rInto \Gamma_5(\Sigma K(\Z/2^{r_1},1)\wedge K(\Z/2^{r_2},1))\rOnto \Z/2.
$$
\end{lemma}
\begin{proof}
Let $A=\Z/2^{r_1}\oplus \Z/2^{r_2}$. Let $X=\Sigma K(\Z/2^{r_1},1)\wedge K(\Z/2^{r_2},1)$. Then $X$ is a retract of $\Sigma K(A,1)\wedge K(A,1)$. From Theorem~\ref{theorem3.2}, $\Gamma_4(X)=\Z/2^{r_1}\otimes\Z/2^{r_2}\otimes\Z/2=\Z/2$ is summand of $\pi_4(X)$ and so
$$
\begin{array}{rcl}
\pi_4(X)&\cong &\Gamma_4(X)\oplus H_4(X)\\
&=&\Gamma_4(X)\oplus H_3(K(\Z/2^{r_1},1)\wedge K(\Z/2^{r_2},1))\\
&\cong& \Gamma_4(X)\oplus \Tor(\Z/2^{r_1},\Z/2^{r_2})\\
&\cong&\Z/2\oplus \Z/2^{\min\{r_1,r_2\}}.\\
\end{array}
$$
The assertion follows from Theorem~\ref{theorem2.2}(2).
\end{proof}
There is a canonical choice of skeleton $\sk_n(K(\Z/2^r,1))$ with $$\sk_n(K(\Z/2^r,1)=\sk_{n-1}(K(\Z/2^r,1)\cup e^n.$$ This induces a choice of skeleton
$$
\sk_n(\Sigma K(\Z/2^{r_1},1)\wedge K(\Z/2^{r_2},1))=\Sigma \bigcup_{i+j\leq n} \sk_i(K(\Z/2^{r_1},1))\wedge \sk_j(K(\Z/2^{r_2},1)).
$$

\begin{lemma}\label{lemma4.4}
Let $r_1,r_2\geq 1$. Let $r=\min\{r_1,r_2\}$. Let $X=\Sigma K(\Z/2^{r_1},1)\wedge K(\Z/2^{r_2},1)$. Then
\begin{enumerate}
\item $\sk_4(X)\simeq M(\Z/2^r,3)\vee S^3$.
\item If $r_1=r_2=1$, then $\sk_5(X)\simeq S^5\vee S^5\vee \Sigma \RP^2\wedge \RP^2$.
\item If $\max\{r_1,r_2\}>1$, then $\sk_5(X)\simeq S^5\vee S^5\vee M(\Z/2^r,3)\vee M(\Z/2^r,4).$
\item
$
\Gamma_5(X)\cong \pi_5(\Sigma M(\Z/2^{r_1},1)\wedge M(\Z/2^{r_2},1))=\left\{
\begin{array}{lcl}
\Z/2\oplus\Z/2&\textrm{ if }&r_1=r_2=1,\\
\Z/4\oplus \Z/2&\textrm{ if } &\min\{r_1,r_2\}=1\\
& & \ \textrm{ and } \max\{r_1,r_2\}>1,\\
\Z/2\oplus \Z/2\oplus\Z/2&\textrm{ if }& r_1,r_2>1.\\
\end{array}\right.
$
\end{enumerate}
\end{lemma}
\begin{proof}
We may assume that $r_1\leq r_2$ and so $r=r_1$. Let $x_i$ ($y_i$) be a basis for $H_i(K(\Z/2^{r_1},1);\Z/2)$ ($H_i(K(\Z/2^{r_2},1);\Z/2)$), which represents the $i$-dimensional cell in the space $K(\Z/2^{r_k},1)$. Then
$$
s^{-1}\tilde H_*(\sk_{n+1}(X);\Z/2)
$$
has a basis given by $x_iy_j$ with $i+j\leq n$ and $i,j\geq 1$. In particular, $s^{-1}\tilde H_*(\sk_4(X);\Z/2)$ has a basis $\{x_1y_1,x_1y_2,x_2y_1\}$ with the Bockstein $\beta_{r_1}(x_2y_1)=x_1y_1$. There is (unique up to homotopy) $2$-local $3$-cell complex with this
homological structure which is given by $S^4\vee M(\Z/2^r,3)$. Thus $\sk_4(X)\simeq S^4\vee M(\Z/2^r,3)$, which is assertion (1).

(2) and (3). Observe that $s^{-1}\tilde H_*(\sk_5(X);\Z/2)$ has a basis $\{x_1y_1,x_1y_2,x_2y_1,x_1y_3,x_2y_2,x_3y_3\}$. Let
$$
j\colon \Sigma M(\Z/2^{r_1},1)\wedge M(\Z/2^{r_2},1)\rInto \sk_5(X)
$$
be the canonical inclusion. For $i=1,2$, the composite
$$
S^5\rTo^g \Sigma^2 K(\Z/2^{r_i},1)\cong \Sigma K(\Z/2^{r_i},1)\wedge S^1\rInto \Sigma K(\Z/2^{r_1},1)\wedge K(\Z/2^{r_2},1),
$$
in which $g$ is map that inducing isomorphism on $H_5(\ ;\Z/2)$ as in Lemma~\ref{lemma4.1}, induces a map
$$
f_i\colon S^5\longrightarrow \sk_5(X).
$$
By inspecting homology, the map
$$
(f_1,f_2,j)\colon S^5\vee S^5\vee \Sigma M(\Z/2^{r_1},1)\wedge M(\Z/2^{r_2},1)\longrightarrow \sk_5(X)
$$
induces an isomorphism on mod $2$ homology and so it is a homotopy equivalent localized at $2$. If $\max\{r_1,r_2\}>1$, then
from decomposition~\ref{equation3.7},
$$
\Sigma M(\Z/2^{r_1},1)\wedge M(\Z/2^{r_2},1)\simeq M(\Z/2^r,3)\vee M(\Z/2^r,4)
$$
and so $\sk_5(X)\simeq S^5\vee S^5\vee M(\Z/2^r,3)\vee M(\Z/2^r,4)$ in this case. Thus assertions (2) and (3) follow.

(4). \textbf{Case I.} $\max\{r_1,r_2\}>1$. By the definition of the Whitehead's functor $\Gamma$,
$$
\begin{array}{rcl}
\Gamma_5(X)&=&\mathrm{Im}(\pi_5(\sk_4(X))\to \pi_5(\sk_5(X))\\
&=&\mathrm{Im}(\pi_5(S^4\vee M(\Z/2^r,3))\to \pi_5(S^5\vee S^5\vee M(\Z/2^r,3)\vee M(\Z/2^r,4)))\\
&=&\pi_5(M(\Z/2^r,3)\vee M(\Z/2^r,4))\\
\end{array}
$$
because
$$
M(\Z/2^r,3)\vee M(\Z/2^r,4)))\simeq \Sigma M(\Z/2^{r_1},1)\wedge M(\Z/2^{r_2},1)=(S^4\vee M(\Z/2^r,3))\cup e^5.
$$
Now it suffices to compute
$$
\pi_5(M(\Z/2^r,3)\vee M(\Z/2^r,4))=\pi_5(M(\Z/2^r,3))\oplus\pi_5(M(\Z/2^r,4)).
$$
It is straight forward to see that $\pi_5(M(\Z/2^r,4))=\Z/2$ represented by the composite
$$
S^5\rTo^{\eta} S^4\rInto M(\Z/2^r,4).
$$

If $r=\min\{r_1,r_2\}=1$, then $\pi_5(M(\Z/2^r,3)=\Z/4$ according to~\cite[Proposition 5.1]{Wu1}.

If $r=\min\{r_1,r_2\}>1$, we compute $\pi_5(M\Z/2^r,3)$. Observe that this is in the stable range and so
$$
\pi_5(M(\Z/2^r,3))\cong \pi_5^s(M(\Z/2^r,3)).
$$
Now we are working in the stable homotopy category. Since $\eta\colon S^5\to S^4$ is of order $2$, there is a map
$$
\tilde\eta\colon M(\Z/2,5)\longrightarrow S^4
$$
such that $\tilde\eta|_{S^5}\simeq \eta$. Since the identity map of $M(\Z/2,5)$ is of order 4 (see for instance~\cite[Theorem 4.4]{Toda}),
there is a commutative diagram
\begin{diagram}
S^5&\rInto^j& M(\Z/2,5)& &\\
   &\ldDashto>{\bar\eta}&\dTo>{\tilde\eta}&&\\
M(\Z/2^r,3)&\rTo^{\mathrm{pinch}}& S^4&\rTo^{[2^r]}&S^4,\\
\end{diagram}
where the bottom row is the cofibre sequence. The composite $\bar\eta\colon S^5\to M(\Z/2^r,3)$ represents an element in $\pi^s_5(M(\Z/2^r,3))$ that maps down to $\pi_5^s(S^4)=\Z/2(\eta)$. Since the map $j\colon S^5\to M(\Z/2,5)$ is of order $2$, the composite $\bar\eta\circ j$ is of order $2$. It follows that
\begin{equation}\label{equation4.3}
\pi_5(M(\Z/2^r,3))\cong \pi_5^s(M(\Z/2^r,3))\cong \pi_5^s(S^4)\oplus \pi_5^s(S^3)\cong \Z/2\oplus\Z/2.
\end{equation}

\noindent\textbf{Case II.} $r_1=r_2=1$. In this case, similar to the above arguments,
\begin{equation}\label{equation4.4}
\Gamma_5(X)=\Gamma_5(\Sigma \RP^2\wedge \RP^2)=\pi_5(\Sigma \RP^2\wedge \RP^2).
\end{equation}
We compute this homotopy group. Note that
$$
\Sigma \RP^2\wedge\RP^2=\sk_4(X)\cup e^5=(S^4\vee M(\Z/2,3)\cup e^5.
$$
There is a cofibre sequence
$$
S^4\rTo^{f} S^4\vee M(\Z/2,3)\rTo^g \Sigma \RP^2\wedge \RP^2,
$$
where the composite
$$
S^4\rTo^{f} S^4\vee M(\Z/2,3)\rTo^{\mathrm{proj.}} S^4
$$
is of degree $2$ because
$$
Sq^1_*\colon H_5(\Sigma \RP^2\wedge \RP^2;\Z/2)=\Z/2\longrightarrow H_4(\Sigma \RP^2\wedge\RP^2;\Z/2)=\Z/2\oplus\Z/2
$$
is not zero, and the composite
$$
S^4\rTo^{f} S^4\vee M(\Z/2,3)\rTo^{\mathrm{proj.}} M(\Z/2,3)
$$
is homotopic to the composite
$$S^4\rTo^{\eta}S^3\rInto^j M(\Z/2,3)$$
because $\pi_4(M(\Z/2,3))=\Z/2$ and
$$
Sq^2_*\colon H_5(\Sigma \RP^2\wedge \RP^2;\Z/2)=\Z/2\longrightarrow H_3(\Sigma \RP^2\wedge\RP^2;\Z/2)=\Z/2
$$
is an isomorphism. Since
$$
\Gamma_5(\Sigma \RP^2\wedge \RP^2)=\pi_5(\Sigma \RP^2\wedge \RP^2),
$$
$$
g_*\colon \pi_5(S^4\vee M(\Z/2,3))\longrightarrow \pi_5(\Sigma \RP^2\wedge \RP^2)
$$
is an epimorphism. By applying the Hilton-Milnor Theorem,
$$
\begin{array}{rcl}
\pi_5(S^4\vee M(\Z/2,3))&\cong& \pi_4(\Omega(S^4\vee M(\Z/2,3)))\\
&\cong & \pi_4(\Omega S^4\times \Omega (M(\Z/2,3))\times\Omega\Sigma (\Omega S^4\wedge \Omega M(\Z/2,3)))\\
&\cong&\pi_4(\Omega S^4)\oplus \pi_4(\Omega (M(\Z/2,3)))\\
&\cong&\pi_5(S^4)\oplus \pi_5(M(\Z/2,3))\\
&\cong&\Z/2\oplus \pi_5(M(\Z/2,3)).\\\
\end{array}
$$
From~\cite[Proposition 5.1]{Wu1}, $\pi_5(M(\Z/2,3))=\Z/4$ generated by the homotopy class of any map $\phi\colon S^5\to M(\Z/2,3)$ such that the composite
$$
S^5\to M(\Z/2,3)\to S^4
$$
is homotopic to $\eta$, the generator for $\pi_5(S^4)=\Z/2$, and, for any such a choice of map $\phi$, the element $2[\phi]$ is given by the homotopy class of the composite
$$
S^5\rTo^{\eta}S^4\rTo^{\eta}S^3\rInto^{j} M(\Z/2,3).
$$
From the fact that $g\circ f\simeq\ast$, the composite
$$
\pi_5(S^4)\rTo^{f_*}\pi_5(\pi_5(S^4\vee M(\Z/2,3)))\rTo^{g_*} \pi_5(\Sigma \RP^2\wedge \RP^2)
$$
is zero. Observe that
$$
f_*(\eta)=2\eta+[j\circ\eta\circ\eta]=2[\phi].
$$
Thus $g_*(2[\phi])=0$ and so $\pi_5(\Sigma \RP^2\wedge \RP^2)$ is a quotient group $\Z/2\oplus\Z/2$. On the other hand, from Theorem~\ref{theorem2.2}(2), there is short exact sequence
$$
\Z/2\rInto \Gamma_5(\Sigma \RP^2\wedge\RP^2)=\pi_5(\Sigma\RP^2\wedge\RP^2)\rOnto \Z/2.
$$
It follows that $\pi_5(\Sigma\RP^2\wedge \RP^2)=\Z/2\oplus\Z/2$. The proof is finished.
\end{proof}

Let $\Len^3(2^r)=\sk_3(K(\Z/2^r,1))$ be the $3$-dimensional lens space.
\begin{lemma}\label{lemma4.5}
Let $r_1,r_2\geq1$. Let
$$
\begin{array}{c}
X_1=\Sigma M(\Z/2^{r_1},1)\wedge M(\Z/2^{r_2},1),\\
X_2=\Sigma \Len^3(2^{r_1})\wedge \Len^3(2^{r_2}),\\
X=\Sigma K(\Z/2^{r_1},1)\wedge K(\Z/2^{r_2},1).\\
\end{array}
$$
Then there is a commutative diagram
\begin{diagram}
&&\Gamma_5(X_1)&\rTo^{\cong}&\pi_5(X_1)&&\\
&&\dTo>{\cong}&&\dInto&&\\
&&\Gamma_5(X_2)&\rInto^{\curvearrowleft}&\pi_5(X_2)&\rOnto^{\curvearrowleft}&\Z/2^{r_1}\oplus\Z/2^{r_2}\\
&&\dTo>{\cong}&&\dOnto&&\dOnto\\
H_6(X)&\rTo&\Gamma_5(X)&\rTo&\pi_5(X)&\rOnto&(\Z/2^{\min\{r_1,r_2\}})^{\oplus 2},\\
\end{diagram}
where the rows are exact and the middle a splitting short exact sequence.
\end{lemma}
\begin{proof}
As in the proof in Lemma~\ref{lemma4.4}, $s^{-1}\tilde H_k(X)$ for $k\leq 6$ has a basis
$$
\{x_iy_j\ | \ i+j\leq 6, \ i,j\geq 1\}.
$$
Thus
$$
\sk_4(X)\subseteq X_1\subseteq \sk_5(X)\subseteq X_2\subseteq \sk_7(X)
$$
and so the commutative diagram follows, where
$$
\Gamma_5(X_1)\cong \Gamma_5(X_2)\cong \Gamma_5(X)
$$
are given by Lemma~\ref{lemma4.4}. Since $\sk_5(X)\subseteq X_2$, $\pi_5(X_2)\to \pi_5(X)$ is onto.

Now we show that the middle row in the diagram splits off. By taking the suspension, there is a commutative diagram of short exact sequences
\begin{diagram}
\Gamma_5(X_2)&\rInto& \pi_5(X_2)&\rOnto& H_5(X_2)\\
\dTo>{\cong}&&\dTo&&\dTo>{\cong}\\
\Gamma_6(\Sigma X_2)&\rInto& \pi_6(\Sigma X_2)&\rOnto&H_6(X_2),\\
\end{diagram}
where the left column is an isomorphism because
$$
\Gamma_5(X_2)\cong \pi_5(X_1)\cong \pi_6(\Sigma X_1)\cong \Gamma_6(\Sigma X_2).
$$
Thus
\begin{equation}\label{equation4.5}
\pi_5(X_2)\cong \pi_6(\Sigma X_2)
\end{equation}
by the $5$-Lemma.

From Lemma~\ref{lemma4.1}, there is a map
$$
g\colon S^5\longrightarrow \Sigma^2 K(\Z/2^r,1)
$$
inducing an isomorphism on $H_5(\ ;\Z/2)$. It follows that
\begin{equation}\label{equation4.6}
\Sigma^2\Len^3(2^r)=\Sigma^2\sk_3(K(\Z/2^r,1))\simeq S^5\vee \Sigma^2 M(\Z/2^{r},1)
\end{equation}
and so
\begin{equation}\label{equation4.7}
\begin{array}{rcl}
\Sigma X_2&=&\Sigma^2\Len^3(2^{r_1})\wedge \Len^3(2^{r_2})\\
&\simeq& (S^5\vee \Sigma^2 M(\Z/2^{r_1},1))\wedge \Len^3(2^{r_2})\\
&\simeq& \Sigma^5 \Len^3(2^{r_2})\vee \Sigma^2\Len^3(2^{r_2})\wedge M(\Z/2^{r_1},1)\\
&\simeq& S^8\vee M(\Z/2^{r_2},6)\vee M(\Z/2^{r_1},6)\vee M(\Z/2^{r_2},3)\wedge M(\Z/2^{r_1},1).\\
\end{array}
\end{equation}
Thus
$$
\pi_6(\Sigma X_2)\cong \Z/2^{r_1}\oplus \Z/2^{r_2}\oplus \Gamma_6(\Sigma X_2)
$$
and hence the result.
\end{proof}

\begin{theorem}\label{theorem4.2}
Let $r>1$. Then
$$
\pi_5(\Sigma K(\Z/2^r,1))\cong \pi_5(\Sigma K(\Z/2^r,1)\wedge K(\Z/2^r,1))\cong \Z/2\oplus \Z/2^r\oplus\Z/2^r.
$$
\end{theorem}
\begin{proof}
Let $X=\Sigma K(\Z/2^r,1)\wedge K(\Z/2^r,1)$. By Lemma~\ref{lemma4.5},
$$
\pi_5(X)\cong \Z/2^r\oplus\Z/2^r\oplus \mathrm{Im}(\Gamma_5(X)\to \pi_5(X)).
$$
From Lemma~\ref{lemma4.4},
$$
\Gamma_5(X)=\Z/2^{\oplus 3}.
$$
By Lemma~\ref{lemma4.1}, the composite
$$
\pi_6(X)\to H_6(X)=\Z/2^r\oplus \Z/2^r \to H_6(X;\Z/2)
$$
is zero. Thus
$$
H_6(X)=\Z/2^r\oplus \Z/2^r\longrightarrow \Gamma_5(X)=\Z/2\oplus\Z/2\oplus\Z/2
$$
detects two copies of $\Z/2$-summands in $\Gamma_5(X)$. The proof is finished.
\end{proof}

\subsection{The Group $\pi_5(\Sigma K(\Z/2^{r_1},1)\wedge K(\Z/2^{r_2},1))$ with $r_1<r_2$.}
Our computation is given by analyzing the cell structure. Let $x_i$ be a basis for $\tilde H_i(K(\Z/2^{r_1},\Z/2))$ and let $y_i$ be a basis for $\tilde H_i(K(\Z/2^{r_2};\Z/2))$. Then $$s^{-1}\tilde H_k(\Sigma K(\Z/2^{r_1},1)\wedge K(\Z/2^{r_2},1), \quad k\leq 6,$$ has a basis $\{x_iy_j\ | \ i+j\leq 6\}$. From the assumption that $r_1<r_2$, the Steenrod operation and Bockstein are indicated by the following diagram
\begin{equation}\label{equation4.8}
\begin{diagram}
k=6\quad&\fbox{$x_2y_3$}&                       & \fbox{$x_1y_4$}&     &               && \fbox{$x_3y_2$}   &                   & \fbox{$x_4y_1$}\\
        &\beta_{r_1}\dTo&\ldDashto^{\beta_{r_2}}_0&\dTo>{Sq^2_*}   &     &                &&\beta_{r_2}\dDashto&\ldTo^{\beta_{r_1}}&\dTo>{Sq^2_*}\\
k=5\quad&\fbox{$x_1y_3$}&                       &                &      &\fbox{$x_2y_2$}&&\fbox{$x_3y_1$}    &                   &\\
        &               &                       &                &\ldTo^{\beta_{r_1}}&  &&                   &                   &\\
k=4\quad&               &                       &\fbox{$x_1y_2$} &               &      &&                   &                &\fbox{$x_3y_1$}\\
        &               &                       &                &               &      &&                   &\ldTo^{\beta_{r_1}}&\\
k=3\quad&               &                       &                &               &      && \fbox{$x_1y_1$},&&\\
\end{diagram}
\end{equation}
where the dash arrows mean that the next Bockstein $\beta_{r_2}$, which comes from $H_*(K(\Z/2^{r_2},1))$, does not actually happen in the Bockstein spectral sequence up to this range.

\begin{lemma}\label{lemma4.6}
Let $r_2>r_1\geq1$ and let $X=\Sigma K(\Z/2^{r_1},1)\wedge K(\Z/2^{r_2},1)$. Then the suspension
$$
E\colon \pi_5(X)\longrightarrow \pi_6(\Sigma X)
$$
is an isomorphism.
\end{lemma}
\begin{proof}
From formula~(\ref{equation4.5}) together with the fact that $\pi_{n-1}(\sk_n(Y))\cong \pi_{n-1}(Y)$,
$$
\pi_5(\sk_6(\Sigma \Len^3(2^{r_1})\wedge \Len^3(2^{r_2})))\longrightarrow \pi_6(\sk_7(\Sigma \Len^3(2^{r_1})\wedge \Len^3(2^{r_2}))).
$$
Notice that
$$
\sk_6(X)=\sk_6(\Sigma \Len^3(2^{r_1})\wedge \Len^3(2^{r_2}))\cup e^6\cup e^6
$$
indicated by the elements $x_1y_4$ and $x_4y_1$ in diagram~(\ref{equation4.8}). Then there is a commutative diagram of right exact sequences
\begin{diagram}
\pi_5(S^5\vee S^5)&\rTo^{f_*}&\pi_5(Z)&\rOnto&\pi_5(\sk_6(X))\\
\dTo>{\cong}&&\dTo>{\cong}&&\dTo>{E}\\
\pi_6(S^6\vee S^6)&\rTo^{\Sigma f_*}&\pi_6(\Sigma Z)&\rOnto&\pi_6(\sk_7(\Sigma X)),\\
\end{diagram}
where $Z=\sk_6(\Sigma \Len^3(2^{r_1})\wedge \Len^3(2^{r_2}))$ and $f\colon S^5\vee S^5\to Y$ is the attaching map for $\sk_6(X)$. The assertion follows by the $5$-lemma.
\end{proof}

\begin{theorem}\label{theorem4.3}
Let $r_2>r_1\geq1$. Then
$$
\pi_5(\Sigma K(\Z/2^{r_1},1)\wedge K(\Z/2^{r_2},1))=\left\{
\begin{array}{lcl}
\Z/2\oplus\Z/4&\textrm{ if }& r_1=1,\ r_2=2,\\
\Z/2\oplus \Z/8&\textrm{ if }& r_1=1,\ r_2\geq 3,\\
\Z/2\oplus \Z/2^{r_1}\oplus \Z/2^{r_1+1}&\textrm{ if }& r_2>r_1>1.\\
\end{array}
\right.
$$
\end{theorem}
\begin{proof}
From Lemma~\ref{lemma4.6}, it suffices to compute $\pi_6(\Sigma^2K(\Z/2^{r_1},1)\wedge K(\Z/2^{r_2},1)).$ Let
$X=\sk_7(\Sigma^2K(\Z/2^{r_1},1)\wedge K(\Z/2^{r_2},1)).$ From splitting formula~(\ref{equation4.7}),
$$
\sk_7(\Sigma^2\Len^3(2^{r_1})\wedge \Len^3(2^{r_2}))\simeq M(\Z/2^{r_2},6)\vee M(\Z/2^{r_1},6)\vee M(\Z/2^{r_1},4)\vee M(\Z/2^{r_1},5).
$$
Let $Y=\sk_7(\Sigma^2\Len^3(2^{r_1})\wedge \Len^3(2^{r_2}))$. Then $s^{-2}\tilde H_*(Y;\Z/2)$ has a basis listed in diagram~(\ref{equation4.8}) excluding the elements $x_1y_4$ and $x_4y_1$. Let $P^{n}(2^r)=M(\Z/2^r,n-1)$. The mod homology $\tilde H_*(P^n(2^r);\Z/2)$ has a basis $u_{n-1}^r$ and $v_{n}^r$ with degrees $|u_{n-1}^r|=n-1$, $v_n^r|=n$ and the Bockstein $\beta_r(v_n^r)=u_{n-1}^r$. Since $X=Y\cup e^7\cup e^7$, there is a cofibre sequence
$$
S^6_1\vee S^6_2\rTo^{f} P^7(2^{r_1})\vee P^7(2^{r_2})\vee P^5(2^{r_1})\vee P^6(2^{r_1})\rTo^{g} X\rTo^{q} S^7_1\vee S^7_2,
$$
where
$$
g_*(u_6^{r_1},v_7^{r_2};u_6^{r_2},v_7^{r_2};u_4^{r_1},v_5^{r_1};u_5^{r_1},v_6^{r_1})=s^2(x_1y_3,x_2y_3;x_3y_1,x_3y_2;x_1y_1,x_2y_1;x_1y_2,x_2y_2)
$$
for catching the corresponding elements in $\tilde H_*(X;\Z/2)$. Here the map $f$ is the attaching map with $f|_{S^6_1}, f|_{S^6_2}$ corresponding to the homological classes $s^2(x_4y_1)$ and $s^2(x_1y_4)$, respectively. Namely, the induced boundary map $q\colon X\to S^7_1\vee S^7_2$ has the homological property that $q_*\colon H_7(X;\Z/2)\to H_7(S^7_1\vee S^7_2)$ is given by
$$
q_*(x_2y_3)=q_*(x_3y_2)=0, \ \ q_*(s^2(x_4y_1))=\iota_1,\ \ q_*(s^2(x_1y_4))=\iota_2,
$$
where $\iota_j$ is the basis for $H_7(S^7_j;\Z/2)$. For $j=1,2$, let $X_j$ be the homotopy cofibre of $f_j$. Then there is a commutative diagram
\begin{equation}\label{equation4.9}
\begin{diagram}
S^6_1\vee S^6_2&\rTo^{f}& P^7(2^{r_1})\vee P^7(2^{r_2})\vee P^5(2^{r_1})\vee P^6(2^{r_1})&\rTo^{g}& X&\rTo^{q}& S^7_1\vee S^7_2\\
\uInto&&\uEq&&\uTo>{\theta_j}&&\uInto\\
S^6_j&\rTo^{f|_{S^6_j}}& P^7(2^{r_1})\vee P^7(2^{r_2})\vee P^5(2^{r_1})\vee P^6(2^{r_1})&\rTo^{g_j}& X_j&\rTo^{q_j}& S^7_j.\\
\end{diagram}
\end{equation}
\noindent\textbf{Statement 1.} \textit{$\theta_{j*}\colon \tilde H_*(X_j;\Z/2)\to \tilde H_*(X;\Z/2)$ is a monomorphism. Moreover,
$$\mathrm{Im}(\theta_{1\ast}\colon H_7(X_1;\Z/2)\to H_7(X;\Z/2))$$
has the basis given by $\{s^2(x_2y_3),s^2(x_3y_2),s^2(x_4y_1)\}$ and
$$\mathrm{Im}(\theta_{2\ast}\colon H_7(X_2;\Z/2)\to H_7(X;\Z/2))$$
has the basis given by $\{s^2(x_2y_3),s^2(x_3y_2),s^2(x_1y_4)\}$. Thus a basis for $\tilde H_*(X_j;\Z/2)$ can be listed in diagram~(\ref{equation4.8}) by removing one element.}
The statement follows immediately by applying mod $2$ homology to diagram~(\ref{equation4.9}), where the only simple computation is given by checking the image of $\theta_{j\ast}$.

\bigskip

\noindent\textbf{Statement 2.} \textit{The composite
$$
\phi_j\colon S^6_j\rTo^{f|_{S^6_j}} P^7(2^{r_1})\vee P^7(2^{r_2})\vee P^5(2^{r_1})\vee P^6(2^{r_1})\rTo^{\mathrm{proj.}} P^7(2^{r_1})
$$
is null homotopic for $j=1,2$.}

\bigskip

Consider the commutative diagram of cofibre sequences
\begin{diagram}
S^6_j&\rTo^{f|_{S^6_j}}& P^7(2^{r_1})\vee P^7(2^{r_2})\vee P^5(2^{r_1})\vee P^6(2^{r_1})&\rTo^{g_j}& X_j\\
\dEq&&\dTo^{\mathrm{proj.}}&&\dTo>{\delta}\\
S^6_j&\rTo& P^7(2^{r_1})&\rTo& Z.\\
\end{diagram}
Then $\dim\tilde H_*(Z;\Z/2)=3$ and $\delta_*\colon H_*(X_j;\Z/2)\to H_*(Z;\Z/2)$ is onto. From diagram~(\ref{equation4.8}), the Bockstein
$$
\beta_t\colon H_7(Z;\Z/2)\to H_6(Z;\Z/2)
$$
is $0$ for $t<r_1$ with the first non-trivial Bockstein given by $\beta_{r_1}$ coming from $\beta_{r_1}(x_2y_3)=x_1y_3$ in diagram~(\ref{equation4.8}). Note that $\pi_6(P^7(2^{r_1}))=\Z/2^{r_1}$ generated by the inclusion $\bar \iota\colon S^6\hookrightarrow P^7(2^{r_1})$. Then the homotopy class
$$
[\phi_j]=k\bar \iota
$$
for some $k\in \Z$. If $k\equiv 1\mod{2}$, then $\dim \tilde H_*(Z;\Z/2)=1$ which contradicts to that $\dim \tilde H_*(Z;\Z/2)=3.$ Thus $k$ must be divisible by $2$. Let $k=2^tk'$ with $k'\equiv 1\mod{2}$ for some $t\geq 1$. If $t<r_1$, then there is a nontrivial Bockstein $\beta_t$ on $\tilde H_*(Z;\Z/2)$ which is impossible from the above. Hence $t\geq r_1$ and so $[\phi_j]=0$ in $\pi_6(P^7(2^{r_1}))=\Z/2^{r_1}$. Statement 2 follows.

\bigskip

\noindent\textbf{Statement 3.} \textit{The composite
$$
\psi\colon S^6_1\rTo^{f|_{S^6_1}} P^7(2^{r_1})\vee P^7(2^{r_2})\vee P^5(2^{r_1})\vee P^6(2^{r_1})\rTo^{\mathrm{proj.}} P^6(2^{r_1})
$$
is null homotopic.}

\bigskip

Consider the commutative diagram of cofibre sequences
\begin{diagram}
S^6_1&\rTo^{f|_{S^6_1}}& P^7(2^{r_1})\vee P^7(2^{r_2})\vee P^5(2^{r_1})\vee P^6(2^{r_1})&\rTo^{g_1}& X_1\\
\dEq&&p\dTo^{\mathrm{proj.}}&&\dTo>{\delta}\\
S^6_1&\rTo& P^6(2^{r_1})&\rTo^{g'}& W.\\
\end{diagram}
Then $\dim \tilde H_*(W;\Z/2)=3$ and $\delta_*\colon \tilde H_(X_1;\Z/2)\to\tilde H_*(W;\Z/2)$ is onto. Moreover, $H_7(W;\Z/2)$ has a basis given by $\delta_*(s^2(x_4y_1))$. By Statement 1, a basis for $\tilde H_*(X_1;\Z/2)$ is listed in diagram~(\ref{equation4.8}) by removing $x_1y_4$. The canonical projection
$$
p\colon P^7(2^{r_1})\vee P^7(2^{r_2})\vee P^5(2^{r_1})\vee P^6(2^{r_1})\longrightarrow P^6(2^{r_1})
$$
has the property that $p_*(u_5^{r_1})=u_5^{r_1}$, $p_*(v_6^{r_1})=v_6^{r_1}$ and $p_*(x)=0$ for $x$ to the other elements in the basis for $\tilde H_*(P^7(2^{r_1})\vee P^7(2^{r_2})\vee P^5(2^{r_1})\vee P^6(2^{r_1});\Z/2)$. In particular, $p_*(v_5^{r_1})=0$. Note that
$$
\begin{array}{rcl}
Sq^2_*\delta_*(s^2(x_4y_1))&=&\delta_*(Sq^2_*(s^2(x_4y_1)))\\
&=&\delta_*(s^2(x_2y_1))\\
&=&\delta_*(g_{1\ast}(v_5^{r_1})\\
&=&g'_*\circ p_*(v_5^{r_1})\\
&=&0.\\
\end{array}
$$
If follows that
$$
Sq^2_*\colon H_7(W;\Z/2)\to H_5(W;\Z/2)
$$
is zero. From the exact sequence
$$
\pi_6(S^5)=\Z/2\rTo^{2^{r_1}}_0\pi_6(S^5)=\Z/2\rTo \pi_6(P^6(2^{r_1})\rTo \pi_5(S^5)=\Z\rTo^{2^{r_1}}\Z,
$$
we have
\begin{equation}\label{equation4.10}
\pi_6(P^6(2^{r_1}))=\Z/2
\end{equation}
generated by the composite
$$
\bar \eta\colon S^6\rTo^{\eta}S^5\rInto P^6(2^{r_1}).
$$
Thus the homotopy class $[\psi]=0$ or $\bar\eta$. If $[\psi]=\bar\eta$, then $Sq^2\colon H_7(W;\Z/2)\to H_5(W;\Z/2)$ is not zero, which is impossible from the above. Hence $[\psi]=0$. This finishes the proof for Statement 3.

\bigskip

\noindent\textbf{Statement 4.} \textit{There is a homotopy decomposition
$$
X\simeq P^7(2^{r_1})\vee T_1\vee T_2,
$$
where $\tilde H_*(T_1;\Z/2)$ and $\tilde H_*(T_2;\Z/2)$ have basis listed by the middle and the right modules in diagram~(\ref{equation4.6}), respectively.}

\bigskip

From Statements 2 and 3, the attaching map $f|_{S^6_1}$ maps into the subspace $P^7(2^{r_2})\vee P^5(2^{r_1})$ up to homotopy because, in the range of $\pi_6$, we have
$$
\pi_6( P^7(2^{r_1})\vee P^7(2^{r_2})\vee P^5(2^{r_1})\vee P^6(2^{r_1}))\cong  \pi_6(P^7(2^{r_1}))\oplus \pi_6(P^7(2^{r_2}))\oplus \pi_6(P^5(2^{r_1}))\oplus \pi_6(P^6(2^{r_1})).
$$
Thus there is a homotopy commutative diagram of cofibre sequences
\begin{equation}\label{equation4.11}
\begin{diagram}
S^6_1&\rTo^{f'}& P^7(2^{r_2})\vee P^5(2^{r_1})&\rTo& T_2\\
\dInto&&\dInto&&\dDashto>{i_1}\\
S^6_1\vee S^6_2&\rTo^f& P^7(2^{r_1})\vee P^7(2^{r_2})\vee P^5(2^{r_1})\vee P^6(2^{r_1})&\rTo^{g}& X\\
\dTo>{\mathrm{proj.}}&&\dTo>{\mathrm{proj.}}&&\dDashto>{q_1}\\
S^6_1&\rTo^{f'}& P^7(2^{r_2})\vee P^5(2^{r_1})&\rTo& T_2.\\
\end{diagram}
\end{equation}
From Statement 2, there is a  homotopy commutative diagram of cofibre sequences
\begin{diagram}
S^6_1\vee S^6_2&\rTo^f& P^7(2^{r_1})\vee P^7(2^{r_2})\vee P^5(2^{r_1})\vee P^6(2^{r_1})&\rTo^{g}& X\\
\dTo&&\dTo>{\mathrm{proj.}}&&\dDashto>{q_2}\\
\ast&\rTo& P^7(2^{r_1})&\rEq& P^7(2^{r_1}).\\
\end{diagram}
Now the composite
$$
P^7(2^{r_1})\vee T_2\rTo^{(g|_{P^7(2^{r_1})},i_1)} X\rTo^{(q_2,q_1)} P^7(2^{2^{r_1}})\vee T_2
$$
is a homotopy equivalence by inspecting the homology and hence the statement.

\bigskip

\noindent\textit{Computation of the Homotopy Group:} From Statement 4, we have
$$
\pi_6(X)\cong\pi_6(P^7(2^{r_1})\vee T_1\vee T_2)\cong \pi_6(P^7(2^{r_1}))\oplus \pi_6(T_1)\oplus \pi_6(T_2)\cong \Z/2^{r_1}\oplus \pi_6(T_1)\oplus\pi_6(T_2).
$$
For computing $\pi_6(T_1)$, since $T_1=P^6(2^{r_1})\cup e^7$, there is a right exact sequence
$$
\pi_6(S^6)=\Z\rTo \pi_6(P^6(2^{r_1}))=\Z/2\rOnto \pi_6(T_1),
$$
where $\pi_6(P^6(2^{r_1}))=\Z/2$ is given in formula~(\ref{equation4.10}). From diagram~(\ref{equation4.6}),
$$
Sq^2_*\colon H_7(T_1;\Z/2)\longrightarrow H_5(T_1;\Z/2)
$$
is an isomorphism and so the attaching map $S^6\to P^6(2^{r_1})$ of $T_1$ is non-trivial. It follows that $\pi_6(T_1)=0$.

Now we compute $\pi_6(T_2)$. From diagram~\ref{equation4.11}, there is a right exact sequence
$$
\pi_6(S^6)=\Z/\rTo^{f'_*}\pi_6(P^7(2^{r_2})\vee P^5(2^{r_1}))=\pi_6(P^7(2^{r_2}))\oplus \pi_6(P^5(2^{r_1}))\rOnto \pi_6(T_2).
$$
Note that a basis for $\tilde H_*(T_2)$ can be listed in the right module of diagram~(\ref{equation4.6}). The composite
$$
\pi_6(S^6)=\Z\rTo^{f'_*}\pi_6(P^7(2^{r_2}))\oplus \pi_6(P^5(2^{r_1}))\rTo^{\mathrm{proj.}} \pi_6(P^7(2^{r_2}))=\Z/2^{r_2}
$$
is of degree $2^{r_1}$ because of the existence of the Bockstein $\beta_{r_1}$. Moveover the composite
\begin{equation}\label{equation4.12}
S^6\rTo^{f'} P^7(2^{r_2})\vee P^5(2^{r_1})\rTo^{\mathrm{proj.}} P^5(2^{r_1})\rTo^{\mathrm{pinch}} S^5
\end{equation}
is homotopic to $\eta$ because of the existence of the Steenrod operation $Sq^2_*$.

\textbf{Case I.} $r_1=1$. According to~\cite[Proposition 5.1]{Wu1}, $\pi_6(P^5(2))=\Z/4$ generated by the homotopy class of any may $S^6\to P^5(2)$ such that the composite $S^6\to P^5(2)\to S^5$ is $\eta$. It follows that there is a right exact sequence
$$
\Z\rTo^{f'_*=(2^{r_1},\lambda)}\Z/2^{r_2}\oplus \Z/4\rOnto\pi_6(T_2),
$$
where $\lambda\colon \Z\to \Z/4$ is an epimorphism. Thus
\begin{equation}\label{equation4.13}
\pi_6(T_2)=\left\{
\begin{array}{rcl}
\Z/4&\textrm{ if }& r_2=2,\\
\Z/8&\textrm{ if }&r_2\geq 3.\\
\end{array}
\right.
\end{equation}

\textbf{Case II.} $r_1>1$. From formula~\ref{equation4.3}, we have $\pi_6(P^5(2^{r_1}))=\Z/2\oplus\Z/2$. Since the composite in~(\ref{equation4.12}) is essential, the composite
$$
\pi_6(S^6)=\Z\rTo^{f'_*}\pi_6(P^7(2^{r_2}))\oplus \pi_6(P^5(2^{r_1}))\rTo^{\mathrm{proj.}} \pi_6(P^5(2^{r_1}))=\Z/2\oplus\Z/2
$$
is nontrivial and so there is right exact sequence
$$
\Z\rTo^{(2^{r_1},\lambda)}\Z/2^{r_2}\oplus (\Z/2\oplus\Z/2)\rOnto \pi_6(T_2)
$$
with $\lambda\colon \Z\to \Z/2\oplus\Z/2$ nontrivial. It follows that
\begin{equation}\label{equation4.14}
\pi_6(T_2)=\Z/2^{r_1+1}\oplus\Z/2 \textrm{ for } r_1>1.
\end{equation}
 The proof is finished now.
\end{proof}

\subsection{The Group $\pi_5(\Sigma^2K(\Z/2^r,1))$.} We use the spectral sequence induced from Carlsson's construction for computing this group.
Let $A$ be an abelian group and
$$ 0\to A_1\buildrel{\delta}\over\to A_0\to A\to 0
$$
a two-step flat resolution of $A$, i.e. $A_0$ is a free abelian
group. The diagram (\ref{bldia}) implies that there is a natural
isomorphism
$$
\pi_4(\Sigma^2K(A,1))\simeq A\widetilde \otimes A,
$$
where $\widetilde \otimes^2$:
$$\widetilde\otimes^2(A)=
A\widetilde \otimes A:=A\otimes A/(a\otimes b+b\otimes a,\ a,b\in
A).
$$
Given a free abelian group $\bar A,$ theorem \ref{theorem2.2} (2)
implies the following natural exact sequence:
$$
\xyma{\Gamma_5(\Sigma^2K(\bar A,1)) \ar@{^{(}->}[r] \ar@{=}[d]&
\pi_5(\Sigma^2K(\bar A,1))\ar@{->>}[r] \ar@{=}[d]& H_5(\Sigma^2K(\bar A,1))\ar@{=}[d]\\
\bar A\widetilde\otimes \bar A\otimes \mathbb Z/2\oplus
\Lambda^2(\bar A)\ar@{^{(}->}[r] & \pi_5(\Sigma^2K(\bar
A,1))\ar@{->>}[r] & \Lambda^3(\bar A)}
$$

The spectral sequence (\ref{specseq}) for $n=2$, gives the
following diagram of exact sequences:
\begin{equation}\label{diagram5}
\xyma{ L_1\Lambda^3(A)\ar@{->}[d] \\ A\widetilde\otimes A\otimes
\mathbb Z/2\oplus \Lambda^2(A) \ar@{^{(}->}[r] \ar@{->}[d] &
\pi_5(\Sigma^2K(A,1)) \ar@{->>}[r] &
L_1\widetilde\otimes^2(A) \ar@{=}[d]\\
\pi_0(\pi_5(\Sigma^2K(N^{-1}(A_1\buildrel\delta\over\to
A_0),1)))\ar@{^{(}->}[r] \ar@{->>}[d] &
\pi_5(\Sigma^2K(A,1))\ar@{->>}[r] \ar@{=}[u] &
\pi_1(\pi_4\Sigma^2K(N^{-1}(A_1\buildrel\delta\over\to A_0),1))\\
\Lambda^3(A)}
\end{equation}
Consider the first derived functor of the functor $\widetilde
\otimes^2$. The short exact sequence
$$
LSP^2(A)\to L\otimes^2(A)\to L\widetilde\otimes^2(A)
$$
in the derived category has the following model:
$$
\xyma{\Lambda^2(A_1) \ar@{^{(}->}[r]^{\delta_2} \ar@{^{(}->}[d] &
A_1\otimes A_0
\ar@{->}[r]^{\delta_1} \ar@{^{(}->}[d] & SP^2(A_0)\ar@{^{(}->}[d] \\
A_1\otimes A_1 \ar@{^{(}->}[r]^{\delta_2'\ \ \ \ \ \ \ \ }
\ar@{->>}[d] & (A_1\otimes A_0)\oplus (A_0\otimes A_1)
\ar@{->>}[d]
\ar@{->}[r]^{\ \ \ \ \ \ \ \ \delta_1'} & A_0\otimes A_0\ar@{->>}[d] \\
SP^2(A_1) \ar@{^{(}->}[r]^{\delta_2''} & A_1\otimes A_0
\ar@{->}[r]^{\delta_1''} & A_0\widetilde\otimes A_0 }
$$
with
\begin{align*}
& \delta_2(a_1\wedge a_1')=a_1\otimes \delta(a_1')-a_1'\otimes \delta(a_1)\\
& \delta_1(a_1\otimes a_0)=a_0\delta(a_1)\\
& \delta_2'(a_1\otimes a_1')=(a_1\otimes \delta(a_1'),-a_1'\otimes
\delta(a_1))\\
& \delta_1'(a_1\otimes a_0,a_1'\otimes a_0')=\delta(a_1)\otimes
a_0+\delta(a_1')\otimes a_0'\\
& \delta_2''(a_1a_1')=a_1\otimes \delta(a_1')+a_1'\otimes
\delta(a_1)\\
& \delta_1''(a_1\otimes a_0)=\partial(a_1)\widetilde\otimes a_0
\end{align*}
for $a_0,a_0'\in A_0,\ a_1,a_1'\in A_1.$ For $n\geq 2,$ looking at
the resolution $\mathbb Z\buildrel{n}\over\to \mathbb Z$ of the
cyclic group $\mathbb Z/n$, we obtain the following representative
of the element $L\widetilde \otimes^2(\mathbb Z/n)$ in the derived
category:
$$
\mathbb Z\buildrel{2n}\over\to \mathbb Z\buildrel{n}\over\to
\mathbb Z/2
$$
In particular,
\begin{align}
& L_1\widetilde\otimes^2(\mathbb Z/2^k)=\mathbb Z/2^{k+1},\ k\geq
1.\label{comp2}
\end{align}
Here $L_1\widetilde\otimes^2$ denotes the first derived functor of
$\widetilde\otimes^2$ (see \ref{subsection2.2}).

We will use the following:
\begin{lemma}\label{wulemma} {(Lemma 2.1 from \cite{Wu})} Let $G_*$ be a
simplicial group and let $n\geq 0$. Suppose that $\pi_0(G_*)$ acts
trivially on $\pi_n(G_*)$. Then the homotopy group $\pi_n(G_*)$ is
contained in the center of $G_n/\mathcal BG_n$, where $\mathcal
BG_n$ is the $n$th simplicial boundary subgroup of $G_n$. \hfill
$\Box$
\end{lemma}

\begin{theorem}\label{theorem4.4}
The homotopy group
$$
\pi_5(\Sigma^2 K(\Z/2^r,1))=\left\{
\begin{array}{lcl}
\Z/8&\textrm{ if }&r=1\\
\Z/2^{r+1}\oplus \Z/2&\textrm{ if }&r>1.\\
\end{array}
\right.
$$
\end{theorem}
\begin{proof} {\bf Case 1: $r=1$.}
The natural epimorphism $\mathbb Z\to \mathbb Z/2$ induces the
homomorphisms $$\pi_n(S^3)=\pi_n(\Sigma^2K(\mathbb Z,1))\to
\pi_n(\Sigma^2K(\mathbb Z/2,1))=\pi_n(\Sigma^2\mathbb RP^\infty),\
n\geq 1.$$ The diagram (\ref{diagram5}) together with
(\ref{comp2}) implies the following short exact sequences:
\begin{equation}\label{rpi}
\xyma{\mathbb Z/2 \ar@{=}[r] \ar@{=}[d] & \pi_5(S^3)\ar@{->}[d]\\
\mathbb Z/2\ar@{^{(}->}[r]& \pi_5(\Sigma^2\mathbb RP^\infty)
\ar@{->>}[r] & \mathbb Z/4}
\end{equation}
Consider this map simplicially, at the level of the natural map between the
Carlsson constructions $F(S^2)=F^{\mathbb Z}(S^2)\to F^{\mathbb
Z/2}(S^2):$
\begin{center}
\begin{tabular}{ccccccccccc}
$F(S^2)_4$ & $\begin{matrix}\longrightarrow\\[-3.5mm]\ldots\\[-2.5mm]\longrightarrow\\[-3.5mm]\longleftarrow\\[-3.5mm]
\ldots\\[-2.5mm]\longleftarrow
\end{matrix}$ & $F(S^2)_3$ & $\begin{matrix}\longrightarrow\\[-3.5mm]\longrightarrow\\[-3.5mm]\longrightarrow\\[-3.5mm]\longrightarrow\\[-3.5mm]\longleftarrow\\[-3.5mm]
\longleftarrow\\[-3.5mm]\longleftarrow
\end{matrix}$ & $\mathbb Z$ \\
$\downarrow$ & & $\downarrow$ & & $\downarrow$ \\
$F^{\mathbb Z/2}(S^2)_4$ & $\begin{matrix}\longrightarrow\\[-3.5mm]\ldots\\[-2.5mm]\longrightarrow\\[-3.5mm]\longleftarrow\\[-3.5mm]
\ldots\\[-2.5mm]\longleftarrow
\end{matrix}$ & $F^{\mathbb Z/2}(S^2)_3$ & $\begin{matrix}\longrightarrow\\[-3.5mm]\longrightarrow\\[-3.5mm]\longrightarrow\\[-3.5mm]\longrightarrow\\[-3.5mm]\longleftarrow\\[-3.5mm]
\longleftarrow\\[-3.5mm]\longleftarrow
\end{matrix}$ & $F^{\mathbb Z/2}(S^2)_2$
\end{tabular}
\end{center}
Here $F^{\mathbb Z/2}(S^2)_k$ is the free product of
$\binom{k}{2}$ copies of $\mathbb Z/2$. In particular $$F^{\mathbb
Z/2}(S^2)_4=\langle s_js_i(\sigma)\ 0\leq i<j\leq 3|\
(s_js_i(\sigma))^2=1\rangle$$

Using the description of the element (\ref{susp5el}), we see that
the simplicial cycle which defines the image of $\pi_5(S^3)$ in
$\pi_5(\Sigma^2\mathbb RP^\infty)$ can be chosen of the form
$$
[[s_2s_1(\sigma),s_1s_0(\sigma)],[s_2s_1(\sigma),s_2s_0(\sigma)]]\in
F^{\mathbb Z/2}(S^2)_4
$$
With the help of lemma \ref{wulemma}, we have
\begin{multline*}
[[s_2s_1(\sigma),s_1s_0(\sigma)],[s_2s_1(\sigma),s_2s_0(\sigma)]]=[[(s_2s_1(\sigma),s_1s_0(\sigma)],(s_2s_1(\sigma)s_2s_0(\sigma))^2]\equiv\\
[[s_2s_1(\sigma),s_1s_0(\sigma)],(s_2s_1(\sigma),s_2s_0(\sigma))]^2\mod
\mathcal B^{\mathbb Z/2}(S^2)_4
\end{multline*}
since
$[[s_2s_1(\sigma),s_1s_0(\sigma)],s_2s_1(\sigma)s_2s_0(\sigma)]$
is a cycle in $F^{\mathbb Z/2}(S^2).$ That is, the image of the
element $\pi_5(S^3)$ is divisible by 2 in $\pi_5(\Sigma^2\mathbb
RP^\infty)$. The diagram (\ref{rpi}) implies the result.
\vspace{.5cm}

\noindent{\bf Case 2: $r>1$.} Now the diagram (\ref{diagram5})
together with (\ref{comp2}) implies the following short exact
sequence
\begin{equation}\label{aaa2}
0\to \mathbb Z/2\to \pi_5(\Sigma^2 K(\Z/2^r,1))\to \mathbb
Z/2^{r+1}\to 0 \end{equation} Therefore, $\pi_5(\Sigma^2
K(\Z/2^r,1))$ is either $\mathbb Z/2^{r+2}$ or $\mathbb
Z/2^{r+1}\oplus \Z/2$.

By theorem \ref{theorem2.2} (2), the Whitehead exact sequence for
$\Sigma^2K(A,1)$ has the following form:
\begin{equation}\label{newdi1}
\xyma{& A\widetilde\otimes A\otimes \Z/2\oplus
\Lambda^2(A)\ar@{^{(}->}[d]\\
H_4(A) \ar@{->}[r] & \Gamma_5(\Sigma^2K(A,1)) \ar@{->>}[d]
\ar@{->}[r] &
\pi_5(\Sigma^2K(A,1))\ar@{->}[r] & H_3(A)\\
& Tor(A,\Z/2) } \end{equation}
For $A=\Z/2^r$ it is of the
following form:
\begin{equation}\label{aaa1}
\xyma{\Z/2 \ar@{^{(}->}[d]\\
\Gamma_5(\Sigma^2K(\Z/2^r,1))\ar@{^{(}->}[r] \ar@{->>}[d] &
\pi_5(\Sigma^2
K(\Z/2^r,1)) \ar@{->>}[r] & \mathbb Z/2^r\\
 \mathbb Z/2}
\end{equation}
The natural projection $\mathbb Z/2^r\twoheadrightarrow \mathbb
Z/2$ induces the map \begin{equation}\label{aaa3} \xyma{\mathbb
Z/2 \ar@{^{(}->}[d]
\ar@{->}[r]^{\simeq}& \Z/2 \ar@{^{(}->}[d]\\
\Gamma_5(\Sigma^2K(\Z/2^r,1)) \ar@{->>}[d]
\ar@{->}[r] & \Gamma_5(\Sigma^2 \mathbb RP^\infty)\ar@{->>}[d]\\
\mathbb Z/2 \ar@{->}[r]^0 & \mathbb Z/2} \end{equation} where the
lower map is zero since the induced map $Tor(\mathbb Z/2^r,\mathbb
Z/2)\to Tor(\mathbb Z/2,\mathbb Z/2)$ is zero. The fact that
$\pi_5(\Sigma^2 \mathbb RP^\infty)=\Z/8$ together with diagram
(\ref{newdi1}) implies that $\Gamma_5(\Sigma^2 \mathbb
RP^\infty)=\mathbb Z/4$. Hence
$\Gamma_5(\Sigma^2K(\Z/2^r,1))=\Z/2\oplus \Z/2,$ since there is no
any endomorphism $\mathbb Z/4\to \mathbb Z/4$ with zero map on
quotients $\mathbb Z/2\to \mathbb Z/2$ (as in diagram
(\ref{aaa3})). The diagram (\ref{aaa1}) and exact sequence
(\ref{aaa2}) implies that $\pi_5(\Sigma^2K(\mathbb
Z/2^r,1))=\mathbb Z/2^{r+1}\oplus \Z/2$.
\end{proof}

\subsection{Applications}

\begin{prop}\label{proposition4.7}
The group $\Gamma_5(\Sigma \mathbb RP^\infty)=\Z/2\oplus\Z/2\oplus\Z/2$.
\end{prop}
\begin{proof}
Consider the Whitehead exact sequence
\begin{multline*}
\pi_6(\Sigma \RP^\infty)\rTo^{h}_{0} H_6(\Sigma
\RP^\infty)=\Z/2\rTo \Gamma_5(\Sigma\RP^\infty)\rTo\\ \pi_5(\Sigma
\RP^\infty)=\Z/2\oplus\Z/2\rTo H_5(\Sigma \RP^\infty)=0,
\end{multline*}
where $\pi_5(\Sigma\RP^\infty)=\Z/2\oplus\Z/2$ by Theorem~\ref{theorem4.1} and the Hurewicz homomorphism
$$h_6\colon \pi_6(\Sigma\RP^\infty)\to H_6(\Sigma\RP^\infty)$$ is zero because, otherwise, it would induces a splitting of $\Sigma \RP^5$ which is impossible by inspecting the Steenrod operation on mod $2$ homology. Thus the order of $\Gamma_5(\Sigma \RP^\infty)$ is $8$. We have to determine the group $\Gamma_5(\Sigma\RP^\infty)$. By the definition,
$$
\Gamma_5(\Sigma\RP^\infty)=\mathrm{Im}(\pi_5(\Sigma \RP^3)\to \pi_5(\Sigma \RP^4)).
$$
Thus the inclusion $\Sigma \RP^3\to \Sigma \RP^\infty$ induces an epimorphism
$$
\Gamma_5(\Sigma \RP^3)\rOnto \Gamma_5(\Sigma\RP^\infty).
$$
Note that $\RP^3=SO(3)$ and so, by the Hopf fibration,
$$
\begin{array}{rcl}
\pi_5(\Sigma SO(3))&\cong& \pi_5(BSO(3))\oplus \pi_5(\Sigma SO(3)\wedge SO(3))\\
&\cong& \pi_4(SO(3))\oplus \pi_5(\Sigma \RP^3\wedge \RP^3)\\
&\cong&\Z/2\oplus\pi_5(\Sigma\RP^3\wedge \RP^3).\\
\end{array}
$$
From Lemmas~\ref{lemma4.4} and~\ref{lemma4.5},
$$
\begin{array}{rcl}
\pi_5(\Sigma\RP^3\wedge \RP^3)&\cong &\Gamma_5(\Sigma\RP^3\wedge \RP^3)\oplus\Z/2\oplus\Z/2\\
&\cong& \pi_5(\Sigma \RP^2\wedge\RP^2)\oplus\Z/2\oplus\Z/2\\
&\cong&\Z/2\oplus\Z/2\oplus\Z/2\oplus\Z/2.\\
\end{array}
$$
It follows that its quotient $\Gamma_5(\Sigma\RP^\infty)$ must be an elementary $2$-group and so hence the result.
\end{proof}

\begin{prop}\label{proposition4.8}
For the suspended projective spaces,
$$
\pi_5(\Sigma\RP^n)=\left\{
\begin{array}{lcl}
\Z/2&\textrm{ if }& n=1,\\
\Z/2^{\oplus 3}&\textrm{ if }& n=2,\\
\Z/2^{\oplus 5}&\textrm{ if }&n=3,\\
\Z/2^{\oplus 3}&\textrm{ if }&n=4,\\
\Z/2^{\oplus 2}&\textrm{ if }&3\leq n\leq \infty.\\
\end{array}\right.
$$
\end{prop}
\begin{proof}
When $n=1$, $\pi_5(S^2)=\Z/2$ from Toda's table\cite{Toda2}. When $n=2$, $\pi_5(\Sigma \RP^2)=\Z/2^{\oplus 3}$ is given in ~\cite[Theorem 6.36]{Wu1}. When $n=3$, $\pi_5(\Sigma\RP^3)$ has been computed in Proposition~\ref{proposition4.7}. For $n\geq 4$, since $\sk_6(\Sigma\RP^\infty)=\Sigma\RP^5$,
$$
\pi_5(\Sigma \RP^n)\cong\pi_5(\Sigma \RP^\infty)=\Z/2\oplus\Z/2
$$
by Theorem~\ref{theorem4.2}. The remaining case is $\pi_5(\Sigma \RP^4)$.
Let $F$ be the homotopy fibre of the pinch map $\Sigma \RP^6\rTo \Sigma \RP^6/\RP^4=M(\Z/2,6)$. By inspecting the Serre spectral sequence to the fibre sequence
$$
\Omega M(\Z/2,6)\rTo F\rTo \Omega \Sigma \RP^6,
$$
the canonical injection $j\colon \Sigma \RP^4\to F$ induces an isomorphism on $H_k(\ ;\Z/2)$ for $k\leq 6$ and so
$$
j_*\colon \pi_k(\Sigma \RP^4)\longrightarrow \pi_k(F)
$$
is an isomorphism for $k\leq 5$. In particular, $\pi_5(\Sigma \RP^4)\cong \pi_5(F)$. From the exact sequence
$$
\pi_5(\Omega M(\Z/2,6))=\Z/2\rTo \pi_5(F)\rTo\pi_5(\Sigma \RP^6)=\Z/2\oplus\Z/2,
$$
the group $\pi_5(F)$ is of order at most $8$ and so is $\pi_5(\Sigma \RP^4)$. From Proposition~\ref{proposition4.7},
$$
\Gamma_5(\Sigma \RP^\infty)=\mathrm{Im}(\pi_5(\Sigma \RP^3)\to \pi_5(\Sigma \RP^4))=\Z/2^{\oplus 3}.
$$
It follows that $\pi_5(\Sigma\RP^4)=\Z/2^{\oplus 3}$ and hence the result.
\end{proof}

\begin{prop}\label{proposition4.9}
$\pi_5(\Sigma K(\Sigma_3,1))\simeq \mathbb Z/2\oplus \mathbb Z/2.$
\end{prop}
\begin{proof}
This follows from the analysis of the map between the Whitehead exact sequences
(\ref{white}) induced by the natural map $\mathbb
Z/2\hookrightarrow \Sigma_3:$
$$
\xyma{H_5(\mathbb Z/2) \ar@{=}[d] \ar@{->}[r] &\Gamma_5(\Sigma
\mathbb RP^\infty) \ar@{->}[r] \ar@{=}[d]& \pi_5(\Sigma \mathbb
RP^\infty) \ar@{->}[r]\ar@{=}[d] & H_4(\mathbb
Z/2) \ar@{=}[d]\\
H_5(\Sigma_3) \ar@{->}[r] & \Gamma_5(\Sigma K(\Sigma_3,1)) \ar@{->}[r]
& \pi_5(\Sigma K(\Sigma_3,1)) \ar@{->}[r] & H_4(\Sigma_3)}
$$
Here the natural isomorphism $\Gamma_5(\Sigma \mathbb RP^\infty)
\to \Gamma_5(\Sigma K(\Sigma_3,1))$ follows from the diagram
$$
\xyma{L_2\Gamma_2^2(\mathbb Z/4\twoheadrightarrow \mathbb Z/2)
\ar@{->}[d] \ar@{=}[r] &
L_2\Gamma_2^2(\mathbb Z/4\twoheadrightarrow \mathbb Z/2) \ar@{->}[d]\\
\Gamma_2^3(\mathbb Z/4\twoheadrightarrow \mathbb Z/2,\ \mathbb
Z/2\hookrightarrow \mathbb Z/4)\ar@{=}[r] \ar@{->}[d] &
\Gamma_2^3(\mathbb Z/4\twoheadrightarrow \mathbb Z/2,\ \mathbb
Z/2\hookrightarrow \mathbb Z/12)\ar@{->}[d]\\
\Gamma_5(\Sigma \mathbb RP^\infty) \ar@{->}[r] \ar@{->>}[d] &
\Gamma_5(\Sigma K(\Sigma_3,1)) \ar@{->>}[d]\\L_1\Gamma_2^2(\mathbb
Z/4\twoheadrightarrow \mathbb Z/2) \ar@{=}[r] &
L_1\Gamma_2^2(\mathbb Z/4\twoheadrightarrow \mathbb Z/2)}
$$
\end{proof}

\section{Relation to K-theory}

As we mentioned in the introduction, there is a natural relation between
the problem considered and algebraic K-theory. Since the
plus-construction $K(G,1)\to K(G,1)^+$ is a homological
equivalence, there is a natural weak homotopy equivalence
$$
\Sigma K(G,1)\to \Sigma (K(G,1)^+)
$$
This defines the natural suspension map:
$$
\pi_n(K(G,1)^+)\to \pi_{n+1}(\Sigma(K(G,1)^+))=\pi_{n+1}(\Sigma
K(G,1))
$$
for $n\geq 1$.

Given a group $G$ and its maximal perfect normal subgroup $P\lhd
G,$ one has natural isomorphism $\pi_n(K(P,1)^+)\simeq
\pi_n(K(G,1)^+),\ n\geq 2$ since $K(P,1)^+$ is homotopy equivalent
to the universal covering space of $K(G,1)^+$.

For a perfect group $G$, the Whitehead exact sequences form the
following commutative diagram:
\begin{equation}\label{susbb}
\xyma{H_4(G) \ar@{=}[d]\ar@{->}[r] & \Gamma_2(H_2(G))\ar@{->>}[d]
\ar@{->}[r] & \pi_3(K(G,1)^+) \ar@{->>}[d] \ar@{->>}[r] & H_3(G)\ar@{=}[d]\\
H_4(G) \ar@{->}[r] & H_2(G)\otimes \mathbb Z/2 \ar@{->}[r] &
\pi_4(\Sigma K(G,1)) \ar@{->>}[r] & H_3(G)}
\end{equation}
Here we will look at the applications of the following two
classical constructions:\\ \\ 1) Let $R$ be a ring and $G=E(R),$
the group of elementary matrices. The group $E(R)$ is perfect and
the plus-construction $K(E(R),1)^+$ also denoted $\tilde K(R)$,
defines the algebraic K-theory of $R$:
$K_n(R)=\pi_n(K(E(R),1)^+),\ n\geq 2.$\\ \\
2) Let $\Sigma_\infty$ be the infinite permutation groups and
$A_\infty$ is the infinite alternating subgroup. There is the
following description of stable homotopy groups of spheres
\cite{Priddy}:
\begin{equation}\label{shg}
\pi_n^S=\pi_n(K(\Sigma_\infty,1)^+)=\pi_n(K(A_\infty,1)^+),\ n\geq
2. \end{equation}

\subsection{} Let $R$ be a ring. In this case, one has the natural
homomorphisms:
$$
K_n(R)\to \pi_{n+1}(\Sigma K(E(R),1)),\ n\geq 2.
$$
For $n=2,$ clearly one has the natural isomorphism:
\begin{equation}\label{equation6.2}
K_2(R)\simeq H_2(E(R))\simeq \pi_3(\Sigma K(E(R),1)).
\end{equation}
It is shown in \cite{Arl} that the map $\Gamma_2(K_2(R))\to
K_3(R)$ factors as
$$
\Gamma_2(K_2(R))\twoheadrightarrow K_2(R)\otimes K_1(\mathbb
Z)\buildrel{\star}\over\to K_3(R),
$$
where $\star$ is the product in algebraic K-theory: $\star:
K_i(S)\otimes K_j(T)\to K_{i+j}(S\otimes T)$. Hence the diagram
(\ref{susbb}) has the following form:
\begin{equation}\label{kdi}
\xyma{H_4(E(R)) \ar@{=}[d]\ar@{->}[r] &
\Gamma_2(K_2(R))\ar@{->>}[d]
\ar@{->}[r] & K_3(R) \ar@{->>}[d] \ar@{->>}[r] & H_3(E(R))\ar@{=}[d]\\
H_4(E(R)) \ar@{->}[r] & K_2(R)\otimes K_1(\mathbb Z)
\ar@{->}[ur]^{\star} \ar@{->}[r] & \pi_4(\Sigma K(E(R),1))
\ar@{->>}[r] & H_3(E(R))}
\end{equation}
and the natural map
\begin{equation}\label{isoktheory}
K_3(R)\to \pi_4(\Sigma K(E(R),1))
\end{equation}
is an isomorphism.
From equations~(\ref{equation6.2}) and~(\ref{isoktheory}) together with the fact that $SL(\mathbb Z)=E(\mathbb Z)$, we have the following:
\begin{theorem}
The natural homomorphism
$$
K_n(R)\longrightarrow \pi_{n+1}(\Sigma K(E(R),1))
$$
is an isomorphism for $n=2,3$. In particular,
$$
\pi_3(\Sigma K(SL(\mathbb Z),1))\cong K_2(\mathbb Z)\cong \mathbb Z/2 \textrm{ and }
$$
$$
\pi_4(\Sigma K(SL(\mathbb Z),1))\cong K_3(\mathbb Z)\cong \mathbb
Z/48.
$$\hfill $\Box$
\end{theorem}

\begin{remark}{\rm The isomorphism (\ref{isoktheory}) and Carlsson construction $F^{E(R)}(S^1)$ gives a
way, for an element of $K_3(R),$ to associate an element from
$F^{E(R)}(S^1)_3=E(R)*E(R)*E(R)$ (uniquely modulo $\mathcal
BF^{E(R)}(S^1)$):
\begin{equation*}
\xyma{ K_3(R) \ar@{~>}[r] \ar@/^25pt/[rr]^{\simeq}_{\;}="d" &
E(R)*E(R)*E(R) & \frac{\mathcal ZF^{E(R)}(S^1)_3}{\mathcal
BF^{E(R)}(S^1)_3} \ar@{^{(}->}[dl]\\ &
\frac{E(R)*E(R)*E(R)}{\mathcal
BF^{E(R)}(S^1)_3}\ar@{->}[u]^{section}}.
\end{equation*}
It is interesting to represent in this way known elements from
$K_3(R)$ for different rings. For $R=\mathbb Z$, $x\in SL(\mathbb
Z)=E(\mathbb Z),$ denote by $x^{(1)},x^{(2)}, x^{(3)}$ the
correspondent elements in the free cube $SL(\mathbb Z)*SL(\mathbb
Z)*SL(\mathbb Z)$. Take the following commuting elements of
$SL(\mathbb Z)$: $$ u=\left(\begin{matrix} 1 & 0 & 0
\\ 0 & -1 & 0 \\ 0 & 0 & 1\end{matrix}\right),\ \ \ \ \  v=\left(\begin{matrix} 1 & 0 & 0
\\ 0 & 1 & 0 \\ 0 & 0 & -1\end{matrix}\right)
$$
The structure of the element (\ref{susp5el}), diagram  (\ref{kdi})
and well-known facts about structure of $K_2(\mathbb Z)$ imply
that, using the above notations, the element
$$ [[u^{(2)},v^{(3)}],[u^{(1)},
v^{(3)}]]
$$
corresponds to the element of order 2 in $K_3(\mathbb Z)$. It
would be interesting to see an element of $SL(\mathbb
Z)*SL(\mathbb Z)*SL(\mathbb Z)$ which corresponds to the generator
of $K_3(\mathbb Z)=\mathbb Z/48$. \hfill $\Box$}
\end{remark}

Consider the case $R=\mathbb Z$ and $n=5$. In this case,
$E(\mathbb Z)=SL(\mathbb Z)$ and we have the following commutative
diagram with exact horizontal sequences:
$$
\xyma{\mathbb Z\oplus (\mathbb Z/2)^2\ar@{=}[d] \ar@{->>}[r]&
(\mathbb Z/2)^3\ar@{=}[d] \ar@{->}[r] & 0 \ar@{=}[d] \ar@{->}[r] &
\mathbb Z/2\ar@{=}[d]\\ H_5SL(\mathbb Z)\ar@{->}[r] \ar@{=}[d] &
\Gamma_4(\tilde K(\mathbb Z))\ar@{->}[r] \ar@{->}[d] & K_4(\mathbb
Z) \ar@{->}[d] \ar@{->}[r] & \mathbb
Z/2 \ar@{=}[d] \ar@{^{(}->}[r] & \mathbb Z/4\ar@{->>}[d]\\
H_5(SL(\mathbb Z))\ar@{->}[r] & \Gamma_5(\Sigma K(SL(\mathbb
Z),1))\ar@{->}[r] & \pi_5(\Sigma K(SL(\mathbb Z),1)) \ar@{->}[r] &
\mathbb Z/2\ar@{->}[r]^0 & \mathbb Z/2 }
$$
and the following commutative diagram:
$$
\xyma{(\mathbb Z/2)^2\ar@{^{(}->}[r]\ar@{=}[d] & (\mathbb
Z/2)^3\ar@{->>}[r] \ar@{=}[d]& \mathbb Z/2\ar@{=}[d]\\
\Gamma_2^2(\Gamma_2(K_2(\mathbb Z))\to K_3(\mathbb
Z))\ar@{^{(}->}[r] \ar@{->}[d] & \Gamma_4(\tilde K(\mathbb Z))
\ar@{->}[d] \ar@{->>}[r] & R_2(K_2(\mathbb Z))\ar@{->}[d]\\
\pi_4(\Sigma K(SL(\mathbb Z),1))\otimes \mathbb Z/2\ar@{^{(}->}[r] &
\Gamma_5(\Sigma K(SL(\mathbb Z),1)) \ar@{->}[r] & Tor(\pi_3(\Sigma
K(SL(\mathbb Z)),1),\mathbb Z/2)\\ \mathbb Z/2 \ar@{=}[u]
\ar@{^{(}->}[r] & (\mathbb Z/2)^2\ar@{=}[u] \ar@{->>}[r] & \mathbb
Z/2 \ar@{=}[u]}
$$
Simple analysis shows that the suspension map $\Gamma_4(\tilde
K(\mathbb Z))\to \Gamma_5(\Sigma K(SL(\mathbb Z),1))$ is an
epimorphism and therefore we have the following theorem:
\begin{theorem}
The Hurewicz homomorphism $$
\pi_5(\Sigma K(SL(\mathbb Z),1))\to H_4(SL(\mathbb Z))=\mathbb Z/2
$$
is an isomorphism.\hfill $\Box$
\end{theorem}

\noindent{\bf Remark.} Since $K_4(\mathbb Z)=0,$ we see that the
natural homomorphism
$$K_4(\mathbb Z)\to \pi_5(\Sigma K(SL(\mathbb Z),1))$$ is not an
isomorphism.

\subsection{} Here we will use (\ref{shg}) for certain computations.
\begin{theorem} Let $A_4$ be the 4-th alternating group. Then
$\pi_4(\Sigma K(A_4,1))=\mathbb Z/4$.
\end{theorem}
\begin{proof}
First recall that\footnote{These computations were done with the
help of HAP-system. The authors thank Graham Ellis for these
computations}
\begin{align*}
& H_1(A_4)=\Z/3,\ \ \ \ \ H_2(A_4)=\Z/2,\ \ \ \ \ \
H_3(A_4)=\Z/6,\ \ \ \ \ H_4(A_4)=0\\
& H_2(A_\infty)=\Z/6,\ \ \ \ \ \ H_3(A_\infty)=\Z/12,\ \ \ \ \ \
\pi_3(\Sigma K(A_4,1))=\Z/6
\end{align*}

Consider the Whitehead exact sequence for the space $\Sigma
K(A_4,1)$:
$${\scriptsize
\xyma{\Gamma_3(\Sigma K(A_4,1))\ar@{^{(}->}[r] & \pi_4(\Sigma
K(A_4,1))\ar@{->}[r] & H_3(A_4)\ar@{->}[r] \ar@{=}[d]&
\Gamma_2(H_1(A_4))\ar@{->}[r] \ar@{=}[d] & \pi_3(\Sigma K(A_4,1))
\ar@{=}[d] \ar@{->>}[r] &
H_2(A_4) \ar@{=}[d]\\
& & \mathbb Z/6\ar@{->}[r] & \mathbb Z/3\ar@{->}[r] & \mathbb Z/6
\ar@{->>}[r] & \mathbb Z/2}}
$$
Since $R_2(\pi_2\Sigma K(A_4,1))=R_2(\mathbb Z/3)=0,$ we have
$\Gamma_3(\Sigma K(A_4,1))=\Gamma_2^2(\mathbb Z/3\hookrightarrow
\mathbb Z/6)$. It follows from the definition of the functor
$\Gamma_2^2$ that it is isomorphic to the pushout
$$
\xyma{\mathbb Z/3\otimes (\mathbb Z/3\oplus \mathbb Z/2)
\ar@{->}[r] \ar@{->}[d] & 0\ar@{->}[d]\\ \mathbb Z/6\otimes
(\mathbb Z/3\oplus \mathbb Z/2) \ar@{->}[r] & \Gamma_2^2(\mathbb
Z/3\hookrightarrow \mathbb Z/6)}
$$
That is, $\Gamma_2^2(\mathbb Z/3\hookrightarrow \mathbb
Z/6)=\mathbb Z/2$ and there is the following short exact sequence:
$$
0\to \mathbb Z/2\to \pi_4(\Sigma K(A_4,1))\to \mathbb Z/2\to 0.
$$
We come to the extension problem: is it $\mathbb Z/2^2$ or
$\mathbb Z/4$?

Consider the monomorphism $A_4\hookrightarrow A_\infty$ and the
map between corresponding Whitehead sequences:
$$
\xyma{\Gamma_3(\Sigma K(A_4,1)) \ar@{^{(}->}[r] \ar@{->}[d]&
\pi_4(\Sigma
K(A_4,1))\ar@{->}[r] \ar@{->}[d]& H_3(A_4) \ar@{->}[r] \ar@{->}[d]& \Gamma_2(H_1(A_4)) \ar@{->}[d]\\
\Gamma_3(\Sigma K(A_\infty,1)) \ar@{->}[r] & \pi_4(\Sigma
K(A_\infty,1))\ar@{->}[r] & H_3(A_\infty) \ar@{->}[r] &
\Gamma_2(H_1(A_\infty))}
$$
which is
\begin{equation}
\xyma{ \mathbb Z/2 \ar@{^{(}->}[r] \ar@{->}[d]& \pi_4(\Sigma
K(A_4,1))\ar@{->}[r] \ar@{->}[d]& \mathbb Z/6 \ar@{->}[r] \ar@{->}[d]& \mathbb Z/6 \ar@{->}[d]\\
\Gamma_2^2(0\to \mathbb Z/6) \ar@{->}[r] & \pi_4(\Sigma
K(A_\infty,1))\ar@{->}[r] & \mathbb Z/12 \ar@{->}[r] & 0}
\end{equation}

It is easy to see that $\Gamma_2^2(0\to \mathbb Z/6)=\mathbb Z/2$
and that $\Gamma_2^2(\mathbb Z/3\hookrightarrow \mathbb Z/6)\to
\Gamma_2^2(0\to \mathbb Z/6)$ is an isomorphism. We obtain the
following diagram:
\begin{equation}\label{dia111}
\xyma{\mathbb Z/2 \ar@{^{(}->}[r] \ar@{=}[d] & \pi_4(\Sigma
K(A_4,1)) \ar@{->}[d] \ar@{->>}[r] & \mathbb Z/2\ar@{^{(}->}[d]\\
\mathbb Z/2 \ar@{->}[r] & \pi_4(\Sigma K(A_\infty,1)) \ar@{->>}[r]
& \mathbb Z/12}
\end{equation}

Now we use the isomorphism (\ref{shg}). Consider the suspension:
$$
K(A_\infty,1)^+\to \Omega\Sigma K(A_\infty,1)^+\simeq \Omega\Sigma
K(A_\infty,1)
$$
and the corresponding map between Whitehead sequences:
$$
\xyma{H_4(A_\infty)\ar@{->}[r] \ar@{=}[d]&
\Gamma_2(\pi_2^S)\ar@{->>}[d] \ar@{->}[r] & \pi_3^S\ar@{->>}[r]
\ar@{->}[d]&
H_3(A_\infty) \ar@{=}[d]\\
H_4(A_\infty)\ar@{->}[r] &\pi_2^S\otimes \mathbb Z/2 \ar@{->}[r] &
\pi_4(\Sigma K(A_\infty,1)) \ar@{->}[r] & H_3(A_\infty)}
$$
Since $$\pi_3^S=\mathbb Z/24,\ \pi_4^S=0,$$ we conclude that the
Whitehead sequence for $K(A_\infty,1)^+$ has the following form:
$$
\xyma{H_4(A_\infty) \ar@{^{(}->}[r] \ar@{=}[d] & \Gamma_2(\pi_2^S)
\ar@{->}[r] \ar@{=}[d]
& \pi_3^S \ar@{->>}[r] \ar@{=}[d] & H_3(A_\infty) \ar@{=}[d]\\
\mathbb Z/2 \ar@{^{(}->}[r] & \mathbb Z/4 \ar@{->}[r] & \mathbb
Z/24 \ar@{->>}[r] & \mathbb Z/12}
$$
We conclude that the map
$$
\pi_3^S\to \pi_4(\Sigma K(A_\infty,1))
$$
is an isomorphism and that the map
$$
H_4(A_\infty)\to \Gamma_2^2(0\to \mathbb Z/6)
$$
is the zero map. The diagram (\ref{dia111}) has the following
form: $$\xyma{\mathbb Z/2 \ar@{^{(}->}[r] \ar@{=}[d] &
\pi_4(\Sigma
K(A_4,1)) \ar@{->}[d] \ar@{->>}[r] & \mathbb Z/2\ar@{^{(}->}[d]\\
\mathbb Z/2 \ar@{^{(}->}[r] & \mathbb Z/24 \ar@{->>}[r] & \mathbb
Z/12}$$ The result follows.
\end{proof}

\end{document}